\documentclass[10pt]{amsart}

\usepackage{graphicx}
\usepackage{amssymb,mathrsfs}

\usepackage{enumerate}        

\usepackage[all,dvips]{xy}    
\UseComputerModernTips        

\newcommand{\dbpw}{\mathcal{D}_{b} \left( \mathcal{P}, w \right)}
\newcommand{\kbpw}{\widetilde{K}_{b} \left( \mathcal{P}, w \right)}
\newcommand{\kapw}{\widetilde{K}_{a} \left( \mathcal{P}, w \right)}
\newcommand{\kpw}{\widetilde{K} \left( \mathcal{P}, w \right)}
\newcommand{\dapw}{\mathcal{D}_{a} \left( \mathcal{P}, w \right)}
\newcommand{\dpw}{\mathcal{D} \left( \mathcal{P}, w \right)}

\newcommand{\di}[3]{\mathcal{D}_{#1} \left( \mathcal{#2}, #3 \right)} 
\newcommand{\kt}[3]{\widetilde{K}_{#1} \left( \mathcal{#2}, #3 \right)} 

\newtheorem{theorem}{Theorem}[section]
\newtheorem{lemma}[theorem]{Lemma}

\newtheorem{proposition}[theorem]{Proposition}
\newtheorem{corollary}[theorem]{Corollary}
\newtheorem{conjecture}[theorem]{Conjecture}
\newtheorem{question}[theorem]{Question}
\theoremstyle{definition}
\newtheorem{definition}[theorem]{Definition}
\newtheorem{example}[theorem]{Example}
\newtheorem{remark}[theorem]{Remark}

\begin{document}

\title[Thompson's Group]{The Action of Thompson's Group on a CAT(0) Boundary}
\author[D.Farley]{Daniel Farley}
      \address{Max Planck Institute for Mathematics\\
               53111 Bonn}
      \email{farley@math.uiuc.edu}

\begin{abstract}  
For a given locally finite CAT(0) cubical complex $X$ with base vertex $\ast$, we define the 
\emph{profile} of a given geodesic ray $c$ issuing from $\ast$ to be the collection of all hyperplanes
(in the sense of \cite{Sag}) crossed by $c$.  We give necessary conditions for a collection
of hyperplanes to form the profile of a geodesic ray, and conjecture that these conditions
are also sufficient.

We show that profiles in diagram and picture complexes can be expressed naturally as infinite 
pictures (or diagrams), and use this fact to describe the fixed points at infinity of the actions  by 
Thompson's groups $F$, $T$, and $V$ on their respective  
CAT(0) cubical complexes.  
In particular, the actions of $T$ and $V$ have no global fixed points.
We obtain a partial description of the fixed set of $F$; it consists, at least, of
an arc $c$ of Tits length $\pi /2$, and any other fixed points of $F$ must have one particular profile,
which we describe.  We conjecture that all of the fixed points of $F$ lie on the arc $c$. 

Our results are motivated by the problem of determining whether $F$ is amenable.
\end{abstract}



\maketitle


\section{Introduction}


Thompson's group $F$ is the group of piecewise linear 
homeomorphisms $h: [0,1] \rightarrow [0,1]$
satisfying:
\begin{enumerate}
\item the finitely many points at which $h$ is non-differentiable are all dyadic rational
numbers, and
\item if $h$ is differentiable at $x_0$, then $h'(x_0) \in \{ 2^{i} \mid i \in \mathbb{Z} \}$.
\end{enumerate}
Thompson also described two other groups, $T$ and $V$,
which are (respectively) the groups of  piecewise linear homeomorphisms $h$ of the circle
$[0,1]/(0=1)$ and the right-continuous bijections $h$ of $[0,1)$; in both cases the functions $h$ 
are required to satisfy (1) and (2).
The survey by Cannon, Floyd, and Parry \cite{CFP} is a useful introduction to all of these groups.

We are interested in the following question:
\begin{question} \label{biggie} \cite{Ross}
Is Thompson's group $F$ amenable?
\end{question}
To explain the original interest in \ref{biggie}, we will need a few definitions.
A group $G$ is \emph{elementary amenable} if it is in the smallest class of groups
that is closed under extensions and direct limits, and contains finite and abelian groups.
A group $G$ is \emph{amenable} if there is a measure $\mu : \mathcal{P}(G) \rightarrow [0,1]$
($\mathcal{P}(G)$ is the power set of $G$) 
such that: i)  $\mu$ is finitely additive; ii) $\mu$ is left invariant; and iii) $\mu(G) = 1$.
We let $EG$, $AG$, and $NF$ denote the classes of elementary amenable groups, amenable groups, 
and groups with no free non-abelian subgroups (respectively).  Von Neumann showed that
$EG \subseteq AG \subseteq NF$.  The problem of determining whether these 
inclusions are proper was posed by Day \cite{Day}.   

Brin and Squier \cite{BS} showed that $F \in NF$.
It is proved in \cite{CFP} that $F \not \in EG$.  
Thus, the existence of $F$ implies that at least one of the inclusions 
$EG \subseteq AG \subseteq NF$ is proper for finitely presented groups:  a positive answer to
\ref{biggie} shows that $F \in AG - EG$, and a negative answer shows that $F \in NF - AG$.  Since at
least 1980, when Geoghegan posed Question \ref{biggie}, Thompson's group was expected by many to be 
an example of a finitely presented non-amenable group with no free subgroups.  

We know today that $EG \subsetneq AG \subsetneq NF$, and that the inclusions are proper for both 
finitely generated and finitely presented groups.
Grigorchuk found examples of finitely generated groups in $AG-EG$ \cite{Gri1} and, in 1998,
finitely presented examples as well \cite{Gri2}.  
In 1980, Ol'shanskii \cite{Ol1} found finitely generated
groups in $NF-AG$.  He and Sapir constructed finitely presented groups in $NF-AG$ in 2002 \cite{OlS}.  

Although the original reason to consider \ref{biggie} is thus obsolete, 
the problem of determining whether $F$ is amenable is still of great interest, and
motivates much of the current work about $F$.

Here we attempt to resolve \ref{biggie} negatively using CAT(0) geometry.  Two earlier results are of vital
importance in this.  First, Adams and Ballmann \cite{AB} showed 
that an amenable group $G$ which acts by isometries
on a locally compact CAT(0) space $X$ must either leave a finite-dimensional flat 
invariant or fix a point at infinity. 
Second, \cite{Me1, Me2} showed that Thompson's
groups $F$, $T$, and $V$ act properly, discretely, and by isometries on 
proper CAT(0) cubical complexes $X_{F}$, $X_{T}$, and $X_{V}$.  If $F$ is amenable, it must
therefore either leave a flat invariant or fix a point at infinity in $X_{F}$.  Elementary properties of $F$
(for instance, the fact that $\oplus_{n=1}^{\infty} \mathbb{Z} \subseteq F$ \cite{CFP}) 
imply that $F$ cannot act properly and freely
on a finite-dimensional flat, 
so we would have a proof that $F$ is non-amenable if we showed that $F$ has no global fixed points
in $\partial X_{F}$.  (Note that the groups $T$ and $V$ are known to be non-amenable, since
both are known to contain non-abelian free subgroups.) 

Most of the effort in this paper goes into describing the spaces at infinity of the locally finite complexes
$X_F$, $X_T$, and $X_V$.  (In fact, our methods apply to all diagram groups and picture groups \cite{Me2}.)
At this point, some background
on $X_F$, $X_T$, and $X_V$ is in order; we restrict our remarks to $X_F$ for the sake of simplicity. 
The constructions in \cite{Me1, Me2} come from the theory of diagram groups, which is
due to Guba and Sapir \cite{GS}.  Each vertex of $X_F$ is labelled by a semigroup diagram, which is essentially
a picture demonstrating how to derive an equality $w_1 = w_2$ between words $w_1$, $w_2$ over a semigroup
presentation.  The group $F$ is 
itself a diagram group, so every element in $F$ can be represented by a semigroup
diagram as well.  The action of $F$ on $X_F$ is given by a natural operation on diagrams:  if 
$x \in F$ and $v \in X_{F}^{0}$, then $x \cdot v$ is obtained by stacking the pictures $x$ and $v$, and then 
``reducing dipoles".  (We refer the reader to Section \ref{diagram} for more specifics, 
or to \cite{GS} for a 
complete introduction.)  

Our description of $\partial X_{F}$ is of the same character.
We represent regions of $\partial X_{F}$ as infinite diagrams, which we call
\emph{profiles}.
The action of $F$ on profiles is determined, as before, by stacking diagrams.  As a result, we can
largely reduce the problem of finding fixed points in $\partial X_F$ to a much easier
algebra problem, which can be handled by a case analysis.  Our main theorem is as follows:

\begin{theorem} \label{biggie2}
Thompson's group $F$ fixes an arc in the boundary $\partial X_F$ of Tits length
$\pi /2$.  Any other fixed points on the boundary of $F$ lie in the profile $\Delta_{\infty}$.
Thompson's groups $T$ and $V$ act without global fixed points on their respective boundaries.
\end{theorem}

This unfortunately leaves \ref{biggie} open.  

The problem of finding any remaining fixed points of the action by $F$ appears to be rather delicate.
In Section \ref{infty}, we give some evidence for and against the existence of additional fixed points
in $\partial X_F$.



The paper is organized as follows.  In Section \ref{CAT0}, we collect various facts about CAT(0)
geometry which will be useful in later sections.  In Section \ref{diagram}, we briefly sketch the
definitions of diagram groups and the cubical complexes, called diagram complexes, on which they act.  
In Section \ref{pict},
we describe regions in the space at infinity of diagram complexes using infinite diagrams, and
describe the action of a diagram group on this space at infinity.  Section \ref{main} contains the
main part of the argument, where it is proved that Thompson's group $F$ fixes the profiles 
$\Delta_{L}$, $\Delta_{R}$, $\Delta_{L-R}$, and $\Delta_{\infty}$, each of which can be described 
by an infinite tree. The groups $T$ and $V$ fix only the profile
$\Delta_{\infty}$.  
In Section \ref{end},
we show that the profiles $\Delta_{L-R}$, $\Delta_{L}$, and $\Delta_{R}$ represent (respectively)
the interior of, the ``left'' endpoint of, and the ``right'' endpoint of an arc of length $\pi /2$
in the Tits metric.  We show moreover that all points on this arc are fixed by $F$.
Finally, in Section \ref{infty} we show that the region at infinity which we call 
$\Delta_{\infty}$  contains no fixed points of $T$ or $V$, even though it is fixed as a set.
As a result, one has a proof that $T$ and $V$ fix no point at infinity.  We also discuss the 
problem of determining whether $F$ fixes any points in $\Delta_{\infty}$.  


\section{Background on CAT(0) Spaces} \label{CAT0}


\subsection{Basic Definitions}


We begin by recalling several basic facts about CAT(0) spaces, all of which are taken directly from 
\cite{BH}.

A metric space $X$ is \emph{geodesic} if, for any $x_1 , x_2 \in X$, there is an
isometric embedding $c: [ 0, d(x_1 , x_2 ) ] \rightarrow X$, called a \emph{geodesic},
such that $c(0) = x_1$ and $c( d(x_1 , x_2 )) = x_2$.  We frequently confuse a geodesic
with its image.  A \emph{geodesic triangle} $\Delta(x,y,z)$
consists of three points $x,y,z \in X$ and choices of geodesics $[x,y]$, $[y,z]$, $[x,z]$
connecting them.  Given such a triangle, it is always possible to find points $\overline{x}$,
$\overline{y}$, $\overline{z}$ in two-dimensional Euclidean space $\mathbb{E}^{2}$  
such that
$d_X ( x,y ) = d_{\mathbb{E}^{2}} ( \overline{x} , \overline{y} )$,  
$d_X ( y,z ) = d_{\mathbb{E}^{2}} ( \overline{y} , \overline{z} )$,  
and $d_X ( x,z ) = d_{\mathbb{E}^{2}} ( \overline{x} , \overline{z} )$.  The triangle 
$\overline{\Delta}(\overline{x}, \overline{y}, \overline{z})$
in $\mathbb{E}^{2}$ determined by $\overline{x}$, $\overline{y}$, and $\overline{z}$ is called
a \emph{comparison triangle} for $\Delta$.  There is a map $h: \Delta \rightarrow \overline{\Delta}$ which
sends sides of $\Delta$ isometrically to the corresponding sides of $\overline{\Delta}$.  We say that the
triangle $\Delta$ satisfies the \emph{CAT(0) inequality} if $d_{X}(a,b) \leq d_{\mathbb{E}^{2}}( h(a), h(b))$
whenever $a,b \in \Delta$.  
A geodesic metric space $X$ is \emph{CAT(0)} if all geodesic triangles in $X$ satisfy
the CAT(0) inequality.  CAT(0) spaces are contractible, and \emph{uniquely geodesic}, i.e., given any two points
$x_1$, $x_2$ in a CAT(0) space $X$, there is a unique geodesic connecting $x_1$ to $x_2$.

If $X$ is an arbitrary metric space, and $c: [0,a] \rightarrow X$, $c' : [0,a'] \rightarrow X$
are geodesic segments satisfying $c(0) = c'(0)$, then  we define the \emph{Alexandrov angle} 
$\angle (c, c')$ as follows:
$$ \angle(c,c') := \lim_{\epsilon \rightarrow 0} 
\sup_{0 <  t, t' < \epsilon} \overline{\angle}_{c(0)} (c(t), c'(t')).$$
Here $\overline{\angle}_{c(0)} (c(t), c'(t'))$ is the angle at $\overline{c(0)}$ in the comparison
triangle $\overline{\Delta}$ for $\Delta \left( c(0), c(t), c'(t') \right)$.  Given
three points $x$, $y$, $z$ 
in a CAT(0) space $X$, we let $\angle_{y}(x,z)$ denote the Alexandrov angle
between the (unique) geodesics $[y,x]$ and $[y,z]$.

The CAT(0) inequality can also be expressed in terms of the Alexandrov angle.  If $\Delta$ is a 
geodesic triangle in the metric space
$X$, then $\Delta$ satisfies the CAT(0) inequality
if and only if each Alexandrov angle in $\Delta$ measures less than the corresponding angle in the
comparison triangle $\overline{\Delta}$.  We say that a geodesic metric space $X$ is CAT(0) if every 
geodesic triangle in $X$ satisfies
this version of the CAT(0) inequality.  Bridson and Haefliger \cite{BH} show that this definition 
of CAT(0) spaces is
equivalent to the earlier one.  

A complete CAT(0) space $X$ has a natural space at infinity $\partial X$, which we now define.
Two geodesic rays $c,c' : [0, \infty) \rightarrow X$ are said to be \emph{asymptotic} if there exists a constant
$K$ such that $d(c(t), c'(t)) \leq K$ for all $t \geq 0$.  The set $\partial X$ of \emph{boundary points}
of $X$ (or \emph{points at infinity}) is the set of equivalence classes of geodesic rays, where two
geodesic rays are equivalent if and only if they are asymptotic.  In practice, we will always use
a basepointed version of this construction.  Fix a point $x \in X$.  We define $\partial X$
to be the set of geodesic rays $c: [0, \infty) \rightarrow X$ issuing from $x$, i.e., satisfying $c(0)=x$.  These
two definitions of $\partial X$ are equivalent in a complete CAT(0) space 
by the following proposition:

\begin{proposition} \cite{BH} \label{parallel} 
If $X$ is a complete CAT(0) space and $c : [0, \infty) \rightarrow X$ is a geodesic ray issuing
from $x$, then for every point $x' \in X$ there is a unique geodesic ray $c'$ which issues from
$x'$ and is asymptotic to $c$.
\qed
\end{proposition}  

If a group $G$ acts by isometries on the CAT(0) space $X$, then it is clear that there is an induced
action on $\partial X$, if we regard the latter as the collection of equivalence classes of geodesic rays
in $X$.  If we use the basepointed version of the construction, then the action $\ast$ can be described as follows:
Let $c \in \partial X$; i.e., $c: [0, \infty) \rightarrow X$ is a geodesic ray and $c(0)=x$.  For an isometry
$g \in G$, $g \ast c$ is the unique geodesic ray issuing from $x$ and asymptotic to the left-translate $g \cdot c$ of $c$.
 

\subsection{Convexity in CAT(0) Spaces}


A subset $C$ of a CAT(0) space $X$ is \emph{convex} if, given 
any two points $x_{1}, x_{2} \in C$, the (unique) geodesic segment $[ x_{1}, x_{2} ]$ 
is contained in $C$.  A function 
$f: X \rightarrow \mathbb{R}$ on a geodesic metric space is \emph{convex} if, for any 
geodesic $c: I \rightarrow X$, the composition $f \circ c$ is convex in the ordinary 
sense, i.e., if, for any $t, t' \in I$ and $s \in [0,1]$, 
$$ (f \circ c) \left( (1-s)t + st' \right) \leq (1-s) \left( f \circ c \right)(t) 
+ s \left( f \circ c \right) \left( t' \right).$$

Bridson and Haefliger \cite{BH} show that there is a natural projection $\pi_{C}: X \rightarrow C$
defined whenever $X$ is a complete CAT(0) space and $C$ is a closed convex subspace.  We collect some
basic properties of $\pi_{C}$ here.   

\begin{proposition} \label{projchar}
Let $C$ be a closed, convex subspace of a complete CAT(0) space $X$.
\begin{enumerate}
\item \cite{BH}   
Let 
$d_{C}: X \rightarrow \mathbb{R}$ be defined by the rule
$$d_{C}(x) = \mathrm{inf}_{y \in C} d(x,y).$$   
The function $d_{C}$ is convex.
\item \cite{BH}
For any $x \in X$, there is a unique point $\pi_{C} (x) \in C$ such that 
$d(x, \pi_{C}(x) ) = d_{C}(x)$.  The function $\pi_{C} : X \rightarrow C$
does not increase distances.
\item \cite{BH}
If $x \in X-C$ and $y \in C- \left\{ \pi(x) \right\}$, then $\angle_{\pi(x)}(x,y) \geq \pi /2$. 
\item Fix $x \in X-C$.  If $y \in C$ satisfies $\angle_{y}( x, y' ) \geq \pi /2$ for
all $y' \in C - \{ y \}$, then $y = \pi_{C}(x)$.
\end{enumerate} 
\end{proposition}

\begin{proof}
(4)  Assume that $y$ satisfies the above condition and $y \neq \pi(x)$.  Consider the comparison triangle 
$\overline{\Delta} \left( \overline{x}, \overline{y}, \overline{\pi(x)} \right)$ for 
the geodesic triangle $\Delta ( x, y, \pi(x) )$ in $X$.  Since the comparison triangle is
non-degenerate by our assumptions, at least one of the comparison angles
$\overline{\angle}_{\overline{y}} \left( \overline{x}, \overline{\pi(x)} \right)$,
$\overline{\angle}_{\overline{\pi(x)}} \left( \overline{x}, \overline{y} \right)$ measures less than
$\pi /2$.  Our assumptions and the CAT(0) inequality imply that 
$\overline{\angle}_{\overline{\pi(x)}} \left( \overline{x}, \overline{y} \right) < \pi /2$.  By the
CAT(0) inequality, $\angle_{\pi(x)}(x,y) < \pi /2$ as well.  This violates (3). 
\end{proof}

We say that a geodesic ray $c: [0, \infty) \rightarrow X$ \emph{crosses} a closed, convex subset $C$
if $(\mathrm{Im} \, c) \cap C$ is a non-empty, compact interval and $(\mathrm{Im} \, c) - C$ is disconnected.

\begin{lemma} \label{monotone}
Let $X$ be a CAT(0) space.
\begin{enumerate}
\item
 Suppose that $C$ is a closed convex subset of $X$,   
$c: [0, \infty) \rightarrow X$
is a geodesic ray which crosses $C$, and    
$(\mathrm{Im} \, c) \cap C = c \left( \left[ t_1, t_2 \right] \right)$.  The function 
$d_{C} \circ c$ is strictly monotonically increasing on $[ t_{2}, \infty)$, and 
$\left( d_{C} \circ c \right)(t) \rightarrow \infty$ as $t \rightarrow \infty$.
\item 
If $c, c' : [0, \infty) \rightarrow X$ are two asymptotic geodesic rays in $X$,
then the function $d( c( \underline{~ ~}), c' ( \underline{~ ~}) ): [0,\infty) \rightarrow [0, \infty)$ is
non-increasing.
\end{enumerate} 
\end{lemma}

\begin{proof}
Both parts are standard exercises using basic properties of convex functions.
Part (1) follows from the fact that $d_{C} \circ c : [0, \infty) \rightarrow [0, \infty)$ is convex, 
$\left( d_{C} \circ c \right) (t_{2}) = 0$,
and $\left( d_{C} \circ c \right) (t') >0$ for some $t' > t_{2}$.  Part (2) follows from the
fact that the function $d( c( \underline{~ ~}), c'( \underline{~ ~}) ): [0, \infty) \rightarrow [0, \infty)$ is convex 
and bounded (see \cite{BH}, page 261). 
\end{proof}


\subsection{CAT(0) Cubical Complexes}


We take the following definition of a cubical complex from \cite{BH}:
\begin{definition} (\cite{BH}, pg. 112) A cubical complex $K$ is the quotient of a disjoint union of
cubes $X = \coprod_{\Lambda} I^{n_{\lambda}}$ by an equivalence relation $\sim$.  The restrictions
$p_{\lambda} : I^{n_{\lambda}} \rightarrow K$ of the natural projection $p: X \rightarrow K = X/\sim$
are required to satisfy:
\begin{enumerate}
\item for every $\lambda \in \Lambda$ the map $p_{\lambda}$ is injective;
\item if $p_{\lambda}( I^{n_{\lambda}} ) \cap p_{\lambda'}(I^{n_{\lambda'}}) \neq \emptyset$ then there
is an isometry $h_{ \lambda, \lambda'}$ from a face $T_{\lambda} \subseteq I^{n_{\lambda}}$ onto a face
$T_{\lambda'} \subseteq I^{n_{\lambda'}}$ such that $p_{\lambda}(x) = p_{\lambda'}(x')$ if and only if
$x' = h_{\lambda, \lambda'}(x)$.
\end{enumerate}
\end{definition}   

Let $x$ and $y$ be points in $X$, and let $l(c)$
denote the length of a path $c$.   
$$ d_{\ell}(x,y) = inf \{ l(c) \mid c(0) = x; c(1) = y \}.$$
The function $d_{\ell}: X \times X \rightarrow [0, \infty]$ defines a metric on any 
cubical complex, called the \emph{length metric} \cite{BH}.  A well-known theorem 
due to Gromov \cite{BH} 
says that the length metric on a cubical complex $X$ is  a CAT(0)  
metric if and only if $X$ satisfies the link condition.  
We avoid recounting the precise statement here, but will 
work exclusively with CAT(0) cubical complexes from now on.    

Let $X$ be a complete CAT(0) cubical complex.  Following \cite{Sag}, define a relation 
$\sim$ on edges of $X$, such that $e_1 \sim e_2$ if and only if $e_1$ and $e_2$ are opposite sides
of a square ($2$-cell)  in $X$.  We will sometimes 
call this relation \emph{simple square equivalence}, although it
is not an equivalence relation.  The transitive, reflexive closure of this relation, also
denoted $\sim$, is called \emph{square equivalence}.  It is clear that square equivalence is an 
equivalence relation.

A \emph{combinatorial hyperplane} in $X$ is an equivalence class of edges under $\sim$.  One 
obtains a geometric hyperplane $H$ as follows: let $M_{e}$ be the set of all midpoints of all edges
square equivalent to $e$.  If $C$ is an arbitrary cube, define $H \cap C$ to be the convex hull
of $M_{e} \cap C$ in the cube $C$.  This description determines $H$.  

We now collect some basic properties of CAT(0) cubical complexes.

\begin{theorem} \label{Sageev} \cite{Sag} Let $X$ be a CAT(0) cubical complex. 
\begin{enumerate}
\item If $J$ is a geometric hyperplane in $X$, then $J$ does not intersect itself and partitions
$X$ into two convex components.
\item If $J_1, \ldots, J_k$ are a collection of geometric hyperplanes in $X$ such that
$J_m \cap J_n \neq \emptyset$ for all $m,n$, then $\bigcap J_i \neq \emptyset$.
\item If $x$ and $y$ are vertices in $X$ connected by a geodesic edge-path $p$ of length $n$, then
$p$ crosses $n$ distinct hyperplanes $J_1, \ldots, J_n$, and these are precisely the hyperplanes
which separate $x$ from $y$.  In particular, any other geodesic edge-path $p'$ connecting $x$ to $y$ must
cross precisely the same hyperplanes.   
\item Each geometric hyperplane $J$ is itself a CAT(0) cubical complex.
\qed
\end{enumerate}
\end{theorem}

The following lemma will be useful in Subsection \ref{profiles}.

\begin{lemma} \label{crossing} Let $X$ be a locally finite CAT(0) cubical complex. 
\begin{enumerate}
\item The closed 
$\frac{1}{2}$-neighborhood of a hyperplane $H$ in $X$ factors isometrically as $H \times [0,1]$. 
\item
Let $c: [0, \infty) \rightarrow X$ be a geodesic ray issuing from a vertex $\ast$ of $X$.   
If $c \left( [0, \infty) \right) \cap C \neq \emptyset$ for some open cube $C$ in $X$, then
$c$ crosses every hyperplane passing through $C$.  
\item  
If $H$ is a hyperplane in $X$, $d_{H}$ is non-constant on an open cell $C$ in $X$, and $H \cap C = \emptyset$, 
then there
exists some hyperplane $H_1$ passing through $C$ such that $H_1 \cap H = \emptyset$.
\item Let $H_1$ and $H_2$ be hyperplanes in the CAT(0) cubical complex $X$, and, for $i=1,2$, let $H_{i}^{+}$, $H_{i}^{-}$ 
be the two open, convex components of $X - H_i$.  If the intersections 
$H_{1}^{+} \cap H_{2}^{-}$, $H_{1}^{+} \cap H_{2}^{+}$, $H_{1}^{-} \cap H_{2}^{-}$,
and $H_{1}^{-} \cap H_{2}^{+}$ are all non-empty, then $H_{1} \cap H_{2}$ is also non-empty. 
\end{enumerate}
\end{lemma}

\begin{proof}
\begin{enumerate}
\item This is a consequence of the proof for Theorem 4.10 (page 611) of \cite{Sag}. 
\item Suppose that $c(t_1) \in C$; let $H$ be any hyperplane passing through $C$.  Thus 
$c(t_1) \in H \times (0,1)$; say $c(t_1) \in H \times \{ t' \}$.  Let $s_1$, $s_2$ be arbitrary numbers such
that $s_1 < t' < s_2$ and $\frac{1}{2} \in (s_1 , s_2)$.  The hyperplanes $H \times \{ s_i \} = H_{s_i}$
$(i=1,2)$ separate $X$ into three distinct connected components, one of which is $H \times (s_1 , s_2 )$.  Note
that $H \times [s_1 , s_2 ]$ contains no vertices of $X$, so $c$ must therefore cross either $H_{s_1}$ or
$H_{s_2}$.  If we assume, without loss of generality, that $c$ crosses $H_{s_1}$, then it must be that
$d_{H_{s_1}} (c(t)) \rightarrow \infty$ monotonically on $[ t_1 , \infty )$, for appropriate $t_1$, by 
Lemma \ref{monotone}.  The function
$d_{H_{s_1}}$ is bounded on $H \times [s_1 , s_2 ]$, so the geodesic ray $c$ eventually leaves
$H \times [s_1 , s_2 ]$, and it cannot cross $H_{s_1}$ a second time, due to the monotonicity of
$d_{H_{s_1}} \circ c$ on $[ t_1 , \infty )$.  It follows that $c$ crosses $H_{s_2}$, and thus also
$H_{1/2} = H$.

\item We prove the contrapositive.  Suppose that $H$ is a hyperplane, $C$ is an open cell such that $H \cap C = \emptyset$, and
every hyperplane $H'$ passing through $C$ satisfies $H' \cap H \neq \emptyset$.  We wish to show that $d_{H}$ is constant on
$C$.

Identify $C$ with $(0,1)^{n}$ and fix a factor of $(0,1)^{n}$ (the last one, without loss of generality).  There is a hyperplane
$H'$ such that $H' \cap C = (0,1)^{n-1} \times \left\{ \frac{1}{2} \right\}$.  Let $x_1$, $x_2 \in C$ be two points in $C$ which differ
only in the last coordinate, say $x_1 = ( c_1 , c_2 , \ldots , c_n )$ and $x_2 = ( c_1 , c_2 , c_3 , \ldots , \hat{c}_{n} )$.
We regard these as points in $H'_{c_n} := H' \times \{ c_n \}$ and $H'_{\hat{c}_n} := H' \times \{ \hat{c}_n \}$, respectively.
Let $x = \left( c_1 , c_2 , \ldots , c_{n-1}, \frac{1}{2} \right)$ be in $H'$, which we identify with 
$H' \times \left\{ \frac{1}{2} \right\}$.

We consider the projection $\pi_{H \cap H'} : H' \rightarrow H \cap H'$.  Let us suppose that
$\pi_{H \cap H'} (x) \in C'$ where $C' = (0,1)^{m}$ is an open cube 
(of dimension at least $2$, since $H$ and $H'$ both
pass through $C'$, and $H \neq H'$).  We make the identifications 
$H' \cap C' = \{ (d_1 , \ldots, d_m ) \in C' \mid d_m = 1/2 \}$ and
$H \cap C' = \{ (d_1, \ldots, d_m ) \in C' \mid d_{m-1} = 1/2 \}$.  

Suppose that $\pi_{H \cap H'} (x) = (d'_1 , d'_2, \ldots, d'_{m-2}, 1/2, 1/2 )$.  The requirement that
$\angle_{\pi_{H \cap H'}(x)}(x,y) \geq \pi/2$ for all $y \neq \pi_{H \cap H'}(x)$ 
in $H \cap H'$ guarantees that $[x, \pi_{H \cap H'}(x) ] \cap C' 
\subseteq \{ d'_1 \} \times \ldots \times \{ d'_{m-2} \} \times (0,1) \times \{ 1/2 \}$ (note: the last
coordinate must be $1/2$ since $[x, \pi_{H \cap H'}(x)] \subseteq H'$).  But it follows from this that
$\angle_{\pi_{H \cap H'}(x)}(x,y) \geq \pi/2$ for all $y \neq \pi_{H \cap H'}(x)$ in $H$.  That is:
$\pi_{H \cap H'} (x) = \pi_{H} (x)$
,  by Proposition \ref{projchar} (4), where $\pi_{H \cap H'} : H' \rightarrow H \cap H'$ and $\pi_{H}: X \rightarrow H$
are the projections.   

Therefore the geodesic segment $[x, \pi_{H}(x)]$ ($= [x, \pi_{H \cap H'}(x)]$) is a subset
of $H'$.  Now we consider the geodesic segments $[x, \pi_{H}(x)] \times \{ c_n \} \subseteq H' \times \{ c_n \}$ and
$[x, \pi_{H}(x)] \times \{ \hat{c}_{n} \} \subseteq H' \times \{ \hat{c}_{n} \}$.  These run parallel
to $[x, \pi_{H}(x)]$, and meet $H$ perpendicularly for similar reasons.  It follows that
$[x_1 , \pi_{H}(x_1)] = [x, \pi_{H}(x)] \times \{ c_n \}$ and $[x_2, \pi_{H}(x_2)] = [x, \pi_{H}(x)] \times \{ \hat{c}_{n} \}$.

This implies that $d_{H}(x_1 ) = d_{H}(x) = d_{H} ( x_2 )$, which implies that
the value of $d_{H}$ is independent of the last coordinate.  Now we can argue
coordinate by coordinate to conclude that $d_H$ is constant on $C$.

\item
Assume that the four intersections in the hypothesis are all non-empty.  If we assume also 
that $H_1 \cap H_2 = \emptyset$, 
then it follows that $\left\{ H_{1}^{+} \cup H_{2}^{+}, H_{1}^{-} \cup H_{2}^{-} \right\}$ is an open cover of $X$.  Now each of the
half-spaces $H_{1}^{+}$, $H_{1}^{-}$, $H_{2}^{+}$, and $H_{2}^{-}$ is a convex subspace of a CAT(0) space, and therefore CAT(0) itself.
It follows that each is contractible.  The same reasoning also applies to the four intersections in the hypothesis: each is CAT(0), and
therefore contractible.  

It then follows that each of the sets $X^{+} = H_{1}^{+} \cup H_{2}^{+}$, $X^{-} = H_{1}^{-} \cup H_{2}^{-}$ is simply connected, since
each is the union of two open contractible sets which intersect in an open contractible set.   The intersection $X^{+} \cap X^{-}$ is the 
union of two disjoint open contractible sets:  
$H_{1}^{+} \cap H_{2}^{-}$ and  $H_{2}^{+} \cap H_{1}^{-}$.  Let $c$ be an arc contained in $X^{+}$ connecting $H_{1}^{+} \cap H_{2}^{-}$
to $H_{2}^{+} \cap H_{1}^{-}$ and meeting each in an open segment.  

We apply van Kampen's theorem to the pieces $X^{-} \cup c$ and $X^{+}$.  The first piece $X^{-} \cup c$ satisfies 
$\pi_{1} \left( X^{-} \cup c \right) \cong \mathbb{Z}$,
while the second is simply connected.  The intersection of these two pieces is the simply connected set 
$\left( H_{1}^{+} \cap H_{2}^{-} \right) \cup \left( H_{2}^{+} \cap H_{1}^{-} \right) \cup c$.  It follows that
$\pi_{1} \left( X^{-} \cup X^{+} \right) = \pi_{1} (X)$ is isomorphic to $\mathbb{Z}$.  The space $X$ is CAT(0), however, 
and therefore contractible.
We have a contradiction.
\end{enumerate} 
\end{proof}


\subsection{Profiles of Geodesic Rays in CAT(0) Cubical Complexes} \label{profiles}

Suppose now that $X$ is locally finite, and let $\ast$ be a vertex, which will serve as a basepoint.
If $H$ is any hyperplane in $X$, let $H^{+}$, the \emph{positive half-space} determined by $H$,
be the open complementary component of $X - H$ that doesn't contain $\ast$;  we let $H^{-}$ be the other
open complementary component of $X-H$. 
If $H_1$
and $H_2$ are hyperplanes in $X$, write $H_{1}^{+} \leq H_{2}^{+}$ if
$H_{2}^{+} \subseteq H_{1}^{+}$.  Clearly $\leq$ is a partial order on positive half-spaces.
We also regard $\leq$ as a partial order on hyperplanes, writing $H_1 \leq H_2$ if
$H_{1}^{+} \leq H_{2}^{+}$.  

For any geodesic ray $c$ in $X$ issuing from $\ast$, define $P(c)$, the \emph{profile of $c$},
to be the collection of all positive half-spaces $H^{+}$ such that $H$ is crossed by $c$.

\begin{proposition} \label{fund}
If $c: [0, \infty) \rightarrow X$ is a geodesic ray issuing from $\ast$, then $P(c)$ is non-empty
and satisfies:  
\begin{enumerate}
\item for any finite subset $\left\{ H_{1}^{+}, \ldots,  H_{n}^{+} \right\} \subseteq P(c)$,
$H_{1}^{+} \cap \ldots \cap H_{n}^{+} \neq \emptyset$;
\item the partially ordered set $\left( P(c), \leq \right)$ has no maximal
elements, and
\item if $H_{1}^{+} \in P(c)$ and $H_{2}^{+} \leq H_{1}^{+}$, then
$H_{2}^{+} \in P(c)$.
\end{enumerate}
\end{proposition}

\begin{proof}
If $H_{1}^{+}$, $H_{2}^{+}$, $\ldots$, $H_{n}^{+}$ are in $P(c)$, then there exist real numbers
$t_1, t_2, \ldots, t_n > 0$ such that $c( [ t_i , \infty) ) \subset H_{i}^{+}$, by Lemma \ref{monotone}.
If $t$ is the largest number in $\{ t_1 , \ldots, t_n \}$, then clearly 
$c( [t, \infty)) \subseteq H_{1}^{+} \cap H_{2}^{+} \cap \ldots \cap H_{n}^{+}$.  This proves that property (1) holds
for $P(c)$.  

It is obvious that property (3) is true of $P(c)$.

Suppose that $H^{+} \in P(c)$ is a maximal element, 
and let $c( [ t, \infty )) \subseteq H^{+}$.  Since $d_{H}(c(t')) \rightarrow  \infty$ as $t' \rightarrow \infty$ by Lemma
\ref{monotone},
we can choose $t$ so that $d_{H}(c(t)) > 1/2$.  We consider the collection $\mathcal{C}$ of all open cells $C$ such that
$c( [t, \infty )) \cap C \neq \emptyset$.  
Note that if $\widehat{C} \cap H \neq \emptyset$,
then any point in $c([ t, \infty)) \cap \widehat{C}$ is at most $1/2$-distant from $H$, so that every cell $C$ in $\mathcal{C}$
satisfies $C \cap H = \emptyset$.   

Let $C \in \mathcal{C}$.  If $d_{H}$ is non-constant on $C$, then, by Lemma \ref{crossing}(3), there is
a hyperplane $H'$ passing through $C$ such that $H' \cap H = \emptyset$.  We claim that 
this implies $H^{+} < (H')^{+}$.  
If $H^{+} \not < (H')^{+}$, i.e., 
if $(H')^{+} \not \subseteq H^{+}$, then $(H')^{+} \cap H^{-} \neq \emptyset$.  We also know
that $ \ast \in H^{-} \cap (H')^{-}$, $H^{+} \cap (H')^{-} \neq \emptyset$ 
(since $C \subseteq H^{+}$ and $C \cap (H')^{-} \neq \emptyset$),
and $H^{+} \cap (H')^{+} \neq \emptyset$ (since the geodesic ray must cross $H'$ by Lemma \ref{crossing}(2)).  Now it follows
from Lemma \ref{crossing}(4) that $H' \cap H \neq \emptyset$, which is a contradiction.  This proves the claim.  
Now $H^{+} < (H')^{+}$ and $c$ crosses $H'$ by Lemma \ref{crossing}(2), which contradicts the maximality of $H^{+}$.

It follows that
$d_{H}$ is constant on all of the cells $C$ in $\mathcal{C}$.   
This contradicts the fact that $d_{H} \circ c$ is strictly
monotonically increasing on $[t, \infty)$.  It follows that property (2) holds.

It is obvious that $P(c)$ is non-empty.
\end{proof} 

From now on, we call a collection of positive half-spaces $\mathcal{H}$ a \emph{profile}
if it is non-empty and satisfies properties (1)-(3) in Proposition \ref{fund}.  I don't know if
every profile $\mathcal{H}$ in this sense is 
realized by a geodesic ray, i.e., if there is some geodesic ray $c$ issuing
from $\ast$ such that $\mathrm{Im} \, c \, \cap H^{+} \neq \emptyset$ if and only if $H^{+} \in \mathcal{H}$.  
I make
the following conjecture.

\begin{conjecture} \label{conj} Let $X$ be a locally finite CAT(0) cubical complex with base vertex $\ast$.
\begin{enumerate}  
\item Every profile is realized by a geodesic ray issuing from $\ast$.
\item The collection of geodesic rays $c$ having a given, fixed profile $\mathcal{H}$
forms a subset of $\partial X$ of diameter less than or equal to $\pi/2$, where the distance
in question is the angle metric (see \cite{BH}).  
\end{enumerate}
\end{conjecture}

Part (2) of Conjecture \ref{conj} implies, in particular, that each profile represents
a contractible subset of $\partial X$, since $\partial X$ is a CAT(1) space with respect to
the angular metric, and sets of diameter less than $\pi$ are contractible in CAT(1) spaces
(Proposition 1.4(4) in \cite{BH}).  
If Conjecture \ref{conj} is true, then the description of $\partial X$ in terms of profiles 
may therefore give a useful homotopical view of the space at infinity, especially if profiles have a convenient description.  
We will give such a description of profiles 
in diagram complexes later in Section \ref{pict}.  

\begin{example}
We give a quick example to show why the most obvious approach to proving Conjecture \ref{conj} (1) fails.
Suppose that $X$ is a CAT(0) cubical complex with base vertex $v$; let $\mathcal{H} = \{ H_{1}^{+}, \ldots, H_{n}^{+}, 
\ldots \}$
be a profile of $X$.  For $n \in \mathbb{N}$, let $v_{n} \in H_{1}^{+} \cap \ldots \cap H_{n}^{+}$.  One might hope
that the sequence $( v_{n} )$ converges to a point at infinity which realizes the profile $\mathcal{H}$. 

Consider $\mathbb{R}^{2}$ endowed with the usual square complex structure in which integer lattice points
are the vertices.  We let $(0,0)$ be the base vertex.  It is not difficult to see that there are precisely
$8$ profiles; four of these profiles are realized by any geodesic ray issuing from $(0,0)$ and travelling
through one of the four open quadrants, and the other four are realized by geodesic rays travelling along the
coordinate axes.  The conjecture thus clearly holds in this case.  If we try to realize the profile corresponding
to the open quadrant $\mathbb{R}^{2,+} = \{ (x,y) \in \mathbb{R}^{2} \mid x,y > 0 \}$ by the above method, however, then we
find that nothing prevents us from choosing all of our points to be of the form $(x,x^{2})$.  Any such sequence would converge
to a point at infinity having the wrong profile.
\end{example}
 

\section{Diagram Groups} \label{diagram}


\subsection{Basic Definitions}


If $A$ is a set (alphabet), then the \emph{free semigroup on $A$}, denoted $A^{+}$,
is the collection of all positive, non-empty words in $A$, with the operation of
concatenation.  Let $\mathcal{P} = \langle \Sigma \mid \mathcal{R} \rangle$ be a semigroup presentation.
Thus, $\Sigma$ is an alphabet and $\mathcal{R} \subseteq \Sigma^{+} \times \Sigma^{+}$ is a collection
of equalities between elements of $\Sigma^{+}$.  We will follow the convention of \cite{Me1} and impose the additional
assumption that no relation of the form $(w,w)$ occurs in $\mathcal{R}$.

We now define pictures over $\mathcal{P}$.  Begin with a \emph{frame} $\partial \left( [0,1]^{2} \right)$, a finite, possibly empty,
collection $\mathcal{T}$ of \emph{transistors}, each homeomorphic to $[0,1]^{2}$, and a finite, non-empty collection $\mathcal{W}$
of \emph{wires}, each homeomorphic to $[0,1]$.  The frame and transistors all have well-defined top, bottom, left, and right sides,
which are the open sides parallel to the coordinate axes, and do not include the corners.  A wire has well-defined initial and terminal
points (i.e., $0$ and $1$, respectively).  A \emph{picture over $\mathcal{P}$}, denoted $\Delta$, is a quotient of
$$\partial \left( [0,1]^2 \right) \coprod \left( \coprod_{T \in \mathcal{T}} T \right) \coprod 
\left( \coprod_{w \in \mathcal{W}} w \right)$$
(for a choice of sets $\mathcal{T}$ and $\mathcal{W}$) by an equivalence relation $\sim$, together with a labelling function
$\ell : \mathcal{W} \rightarrow \Sigma$, satisfying:
\begin{enumerate}
\item The initial point of any given wire is attached either to the bottom of a transistor, or to the top of the frame.
The terminal point of any given wire is attached either  to the bottom of the frame, or to the top of some transistor.
If $w$ is a wire and $T$ is a transistor, then $w \cap T \subseteq \Delta$  is either empty or a singleton set.
\item Let $T_1$ and $T_2$ be transistors.  Write $T_1 < T_2$ if there is some wire $w$ such that the initial point of $w$
is attached to the bottom of $T_1$ and the terminal point of $w$ is attached to the top of $T_2$.  Let $<$ also denote the
transitive closure of the above relation.  The relation $<$ is required to be a strict partial order.
\item The equivalence classes of $\sim$ are either singleton sets or consist of exactly two points, exactly one of which is an endpoint
of a wire.  In other words, the only identifications in 
$$ \partial \left( [0,1]^{2} \right) \coprod \left( \coprod_{T \in \mathcal{T}} T \right) \coprod 
\left( \coprod_{w \in \mathcal{W}} w \right)$$
are generated by the attaching maps of the wires, and no two wires have points in common.  The endpoints of the wires
are called \emph{contacts}.
\item Suppose the top of the transistor $T$ meets the wires $w_{i_1}$, $w_{i_2}$, $\ldots$, $w_{i_m}$, reading from
the left side of $T$ to the right.  Suppose that the bottom of the transistor $T$ meets the wires $w_{j_1}$, $\ldots$,
$w_{j_n}$, again reading from left to right.  The \emph{top label of $T$}, denoted $L_{T}$, is
$$ \ell \left( w_{i_1} \right) \ell \left( w_{i_2} \right) \ldots \ell \left( w_{i_m} \right);$$
the \emph{bottom label of $T$}, denoted $L_{B}$, is
$$ \ell \left( w_{j_1} \right) \ell \left( w_{j_2} \right) \ldots \ell \left( w_{j_n} \right).$$
We require that $( L_{T} , L_{B} ) \in \mathcal{R}$ or $( L_{B}, L_{T}) \in \mathcal{R}$. 
\end{enumerate}
We can define the top and bottom labels of the frame just as we did for a transistor $T$.  If
the top label of the frame is $w_1$ and the bottom label is $w_2$, then $\Delta$ is a
\emph{$(w_1 , w_2)$-picture over $\mathcal{P}$}.  We say that $\Delta$ is a $(w, \ast )$-picture if
the top label of $\Delta$ is $w$, and the bottom label is arbitrary.

Two pictures $\Delta_{1}$ and $\Delta_{2}$ are \emph{isomorphic}, $\Delta_{1} \equiv \Delta_{2}$, if there is a homeomorphism
between them which matches labels and preserves the top-bottom- and left-right-orientations on the frame and transistors.

Given a $(u,v)$-picture $\Delta_1$ and a $(v,w)$-picture $\Delta_2$, one can define the \emph{concatenation}
$\Delta_{1} \circ \Delta_{2}$, which is the $(u,w)$-picture obtained by identifying the bottom of the frame
for $\Delta_1$ with the top of the frame for $\Delta_2$ by a homeomorphism which matches the endpoints of the wires, and
then removing the line segment corresponding to the bottom of $\Delta_1$ in the quotient, while keeping the wires passing
through this line segment intact.

\begin{figure}[!h] 
\begin{center}
\includegraphics{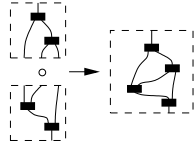}
\end{center}
\caption{On the left, we have two pictures $\Delta_{1}$ and $\Delta_{2}$ (reading from top to bottom); on the right
we have the concatenation $\Delta_{1} \circ \Delta_{2}$.}
\label{concat}
\end{figure}

Figure \ref{concat} illustrates the operation of concatenation in a particular case.  All of the semigroup pictures
in the figure are pictures over the presentation $\langle x \mid x = x^{2} \rangle$.  For this reason, we leave off the labels
of the wires, since the label of each one is $x$.  Note that on the top left is an $(x, x^{3})$-picture, and 
on the bottom left is an $(x^{3}, x)$-picture.  If we denote these pictures $\Delta_{1}$ and $\Delta_{2}$, respectively, then
the $(x,x)$-picture on the right is $\Delta_{1} \circ \Delta_{2}$.

Note that Figure \ref{concat} also illustrates our conventions for drawing pictures in the plane.  If $\Delta$ is a 
semigroup picture, then a function $\rho: 
\Delta \rightarrow \mathbb{R}^{2}$ is a \emph{projection} of $\Delta$ if:
\begin{enumerate}
\item the image of each transistor is a rectangle whose sides are parallel to the coordinate axes.  The map $\rho$
takes the top, left, right, and bottom of any given transistor to the corresponding sides in the image.
\item the image of the frame is an empty rectangle, and the map $\rho$ is again orientation-preserving, in the sense of (1).
The image of $\rho$ is contained inside the image of the frame. 
\item the image of each wire meets any given horizontal line at most once, and
\item $\rho$ is an embedding, except possibly at finitely many double points.  The inverse image of any double point $x$ is
a set of two points on distinct wires $w_1$ and $w_2$.  We assume that the images of $w_1$ and $w_2$ are transverse at $x$.
\end{enumerate}
It is rather clear that all of the defining features of a semigroup picture can be recovered from any of its 
suitably labelled projections.  From
now on, we will usually confuse a picture with any of its projections without further comment.

Two transistors $T_1 < T_2$ form a \emph{dipole} if the top label of $T_1$ is identical (as a word in $\Sigma^{+}$) to the 
bottom label of $T_2$, and the bottom contacts of $T_1$ are paired off by wires in order with the top contacts of $T_2$.  To
\emph{remove a dipole}, delete the transistors $T_1$ and $T_2$ and all wires connecting them, and then glue together in order
the wires that formed top contacts of $T_1$ with those that formed bottom contacts of $T_2$.  The inverse operation is
called \emph{inserting a dipole}.  Two pictures are \emph{equal modulo dipoles}, $\Delta_1 = \Delta_2$, if one can be obtained
from the other by repeatedly inserting and removing dipoles.  A picture is called \emph{reduced} if it contains no dipoles.  Any
equivalence class modulo dipoles contains a unique reduced picture \cite{Me2, GS}.

\begin{figure}[!h] 
\begin{center}
\includegraphics{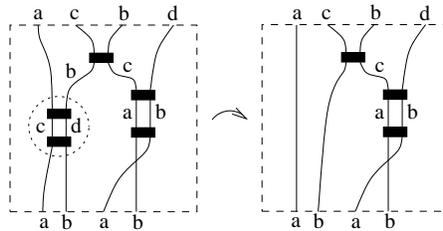}
\end{center}
\caption{The circled transistors in the left half of the figure form a dipole; on the right, we have the
result of removing this dipole.}
\label{dipole}
\end{figure}

In Figure \ref{dipole} we have two $(acbd, abab)$-pictures over the presentation $\mathcal{P} = \langle a, b, c, d \mid
ab=cd, cb=bc, ab=ba \rangle$.  In the left picture, we've circled two transistors which form a dipole.  If we remove this
dipole, we arrive at the picture on the right.  Notice that the two right-most transistors in the right half of the Figure do
not form a dipole:  the top label of the top transistor is $cd$, but the bottom label of the bottom transistor is $ba$. 

For a fixed word $w \in \Sigma^{+}$, the set of all $(w,w)$-pictures over $\mathcal{P}$, modulo dipoles, forms a group
$\dbpw$ under the operation of concatenation.  We will follow \cite{GS} and call $\dbpw$
the \emph{braided diagram group} over $\mathcal{P}$, based at $w$.  (Warning:  the word ``braided" is rather unfortunate.  In fact, 
as the above definition shows, we don't care about any possible braiding of the wires, since equivalence between pictures doesn't
depend on any embedding into an ambient space.  Moreover, there is now a growing literature 
(see for example \cite{Brin})
on a braided version of Thompson's group $V$, which is something quite different from the older group $V$ which we consider here.
Nevertheless, there seems to be no better term.)  
A picture $\Delta$ is \emph{planar} if there is a projection 
$\rho : \Delta \rightarrow \mathbb{R}^{2}$ which is also an embedding.  
The set of all planar $(w,w)$-pictures over $\mathcal{P}$, modulo dipoles, forms a group
$\dpw$, which we will call the \emph{diagram group} over $\mathcal{P}$, based at $w$.  Annular pictures can be defined as follows.
Suppose that $\Delta$ is a picture, and let $\Delta'$ be the space obtained from $\Delta$ by removing the sides of the frame.  We say that
$\Delta$ is \emph{annular} if there is an orientation-preserving immersion of $\Delta'$ into
$A = \{ (x,y) \in \mathbb{R}^{2} \mid 1 \leq x^2 + y^2 \leq 4 \}$ such that:
i) the top of the frame for $\Delta'$ is wrapped around the circle $x^2 + y^2 =1$ once in the counterclockwise direction.  The initial
and terminal points of the top are both mapped to $(1,0) \in A$;
ii) the bottom of the frame for $\Delta'$ is wrapped around the circle $x^2 + y^2 =4$ once in the counterclockwise direction.  The initial
and terminal points of the bottom are both mapped to $(2,0) \in A$;
iii) the only double points of $\rho$ are $(1,0)$ and $(2,0)$; $\rho$ is an embedding otherwise.
The set of all annular pictures over $\mathcal{P}$ is a group $\dapw$, called the \emph{annular diagram group} over 
$\mathcal{P}$, based at $w$. 

Three groups are of special interest to us.  Let $\mathcal{P} = \langle x \mid x = x^{2} \rangle$.  
The groups $\di{}{\mathcal{P}}{x}$, $\di{a}{\mathcal{P}}{x}$, and $\di{b}{\mathcal{P}}{x}$ are, respectively,
Thompson's groups $F$, $T$, and $V$.  The original observation that $\di{}{\mathcal{P}}{x} \cong F$ was due to Victor Guba;
Guba and Sapir (in \cite{GS}) 
sketched the theory of annular and braided diagram groups expressly for the purpose of bringing their techniques
to bear on the study of $T$ and $V$.  Section 6 of \cite{Me2} describes an isomorphism between the groups $F$, $T$, and $V$, and
the corresponding diagram groups.  


\subsection{Diagram Complexes}


If $G$ is a diagram group of the standard, annular, or braided variety, then
a theorem of \cite{Me2} (see also \cite{Me1}) says that $G$ acts properly by isometries on a CAT(0) cubical
complex.  We briefly describe the construction of the cubical complex in this
subsection.

\begin{figure} [!h] 
\begin{center}
\includegraphics{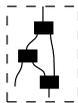}
\end{center}
\caption{This is a vertex in the cubical complex $\kt{b}{\mathcal{P}}{x}$, where
$\mathcal{P} = \langle x \mid x = x^{2} \rangle$.}
\label{vertexpic}
\end{figure}

Fix a braided diagram group $\dbpw$.  We define a complex $\kbpw$, called the
\emph{diagram complex for $\dbpw$}, as follows.  A vertex $v \in \kbpw^{0}$ is an
equivalence class $\sim$ of reduced braided $(w, \ast)$-pictures, where $\Delta_{1} \sim
\Delta_{2}$ if and only if there is some braided permutation picture $\Psi$, such that
$\Delta_{1} \circ \Psi = \Delta_{2}$.  Here a \emph{permutation picture} is one 
with no transistors.  It is convenient to depict a vertex as a $(w,\ast)$-picture
in which all wires 
which would ordinarily be connected to the bottom of the frame have been
cut, as in Figure \ref{vertexpic}.  

\begin{figure} [!h] 
\begin{center}
\includegraphics{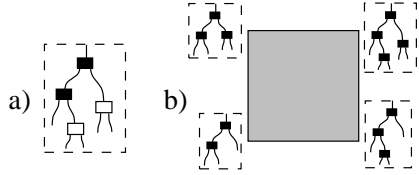}
\end{center}
\caption{a) This notation describes a cube in the complex $\kt{b}{\mathcal{P}}{x}$; b) This is the labelled
cube in $\kt{b}{\mathcal{P}}{x}$ denoted by the picture in a).}
\label{cubepic}
\end{figure}

An $n$-dimensional cube in $\kbpw$ is denoted by a reduced braided $(w,\ast)$-picture $\Delta$ in which
all of the bottom wires have been cut (as above), and $n$ of the maximal
transistors of $\Delta$ have been
drawn as white.  The picture Figure \ref{cubepic}a) denotes a $2$-cube, for example.

If we arbitrarily number the white transistors $1, 2, \ldots, n$, then there is a natural
way to label the vertices of an $n$-cube $[0,1]^{n}$, corresponding to this numbering of
$\Delta$.  Namely, if $(a_1 , \ldots , a_n ) \in \{ 0, 1 \}^{n}$ label $(a_1, \ldots , a_n )$
by the picture $\Delta_{(a_1, \ldots , a_n )}$, where the $i$th transistor is left off if $i=0$
and the $i$th transistor is filled in if $i=1$.  For instance, Figure \ref{cubepic}b) shows how to label the
corners of $[0,1]^{2}$ if $\Delta$ is as in Figure \ref{cubepic}a) and the white transistors are numbered
from left to right.

If we let $\Delta$ vary over all possible isomorphism classes of cube representatives ( where isomorphisms
send white transistors to white transistors), and, for each $\Delta$, choose as above a labelling
of $[0,1]^{n}$ for the appropriate $n$, then $\kbpw$ is the quotient of the resulting labelled cubes by the
equivalence relation which identifies the cubes along faces with the same labels.  It is proved in \cite{Me2} that
$\kbpw$ is a proper CAT(0) cubical complex if $\mathcal{P}$ is a finite presentation, and that
$\dbpw$ acts properly and cellularly on $\kbpw$.  The action is usually not cocompact, and, in particular,
isn't for any 
of the groups $F$, $T$, and $V$. 

We note that entirely similar statements are true for ordinary and annular diagram groups.  It is only
necessary to replace pictures with planar pictures and annular pictures (respectively) in the above
discussion to get the descriptions of $\kpw$ and $\kapw$, respectively.

Lastly, we recall a useful partial order on vertices. 
If $\left[ \Delta_{1} \right]$ and $\left[ \Delta_{2} \right] $ are vertices in a diagram complex, we write
$\left[ \Delta_{1} \right] \leq \left[ \Delta_{2} \right]$ if there exists some picture $\theta$ such that
$\Delta_{1} \circ \theta \equiv \Delta_{2}$.  Note that this means $\Delta_{1} \circ \theta$ and $\Delta_{2}$
are isomorphic before reducing dipoles.  It is not difficult to see that $\leq$ is a well-defined partial order.

Suppose $\mathcal{T}' \subseteq \mathcal{T}_{\Delta}$, where $\mathcal{T}_{\Delta}$ is the collection of transistors
in a picture $\Delta$. We say that $\mathcal{T}'$ is an \emph{initial subset} of $\mathcal{T}_{\Delta}$ if whenever
$T_{1} < T_{2}$ and $T_{2} \in \mathcal{T}'$, then $T_{1} \in \mathcal{T}'$ also.   
We reproduce a lemma from \cite{Me2}.

\begin{lemma} \label{bccc} \cite{Me2}
Let $\Delta$ be a vertex, and let $\mathcal{T}_{\Delta}$ be its set of transistors.  
There is a one-to-one correspondence $\psi$
between initial subsets of $\mathcal{T}_{\Delta}$ and vertices $\Delta_{1}$ satisfying
$\Delta_{1} \leq \Delta$.  The function $\psi$ is order-preserving and has an order-preserving 
inverse, i.e., the initial subsets $\mathcal{T}', \mathcal{T}''$ satisfy $\mathcal{T}' \subseteq \mathcal{T}''$
if and only if $\psi \left( \mathcal{T}' \right) \leq
 \psi \left( \mathcal{T}'' \right)$.
\qed
\end{lemma}

The map $\psi$ in the above lemma is easy to define:  if $\mathcal{T}'$ is an initial subset of transistors, then
$\psi \left( \mathcal{T}' \right)$ is obtained by removing all transistors in $\mathcal{T}_{\Delta} - \mathcal{T}'$,
along with all of their bottom wires.  The result is easily seen to be a vertex.  The argument that the map $\psi$ is 
injective can be extended to prove that the automorphism group of a diagram is trivial, 
at least combinatorially speaking.  That is,
if $\phi : \Delta \rightarrow \Delta$ is a isomorphism, then $\phi$ leaves the frame, each transistor, and each wire invariant,
and restricts to a self-homeomorphism of each of these.  It therefore follows, for instance, that in a concatenation
$\Delta_{1} \circ \Delta_{2}$ of pictures, one can speak of the transistors that were contributed by $\Delta_{i}$ for $i=1,2$, and
this is a well-defined notion even after reducing dipoles.  We shall need this observation in future sections, and use it without
further comment.


\section{Geodesic Profiles in Diagram Complexes} \label{pict}


We now describe profiles in diagram complexes.  Our main goals here are, first,
to describe a profile as an infinite picture of a certain kind, and then to describe
the action on profiles in terms of picture multiplication.

Throughout this section, we use only the complex $\kbpw$, but the discussion carries over
to $\kpw$ and $\kapw$ in an obvious way.

\subsection{Description of Profiles}

The first step is to describe hyperplanes in $\kbpw$.  Recall that a combinatorial hyperplane
is an equivalence class of $1$-cells under the relation $\sim$ of square equivalence.  The square
equivalence relation is generated by simple square equivalence (also denoted $\sim$), 
where two $1$-cells $e_1$, $e_2$ are simple square equivalent  
 if they are opposite faces of a $2$-cell (square).

\begin{figure} [!h] 
\begin{center}
\includegraphics{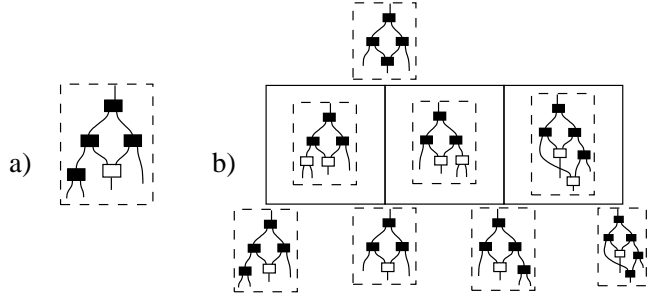}
\end{center}
\caption{a) An edge $\Delta$ in $\kt{b}{\mathcal{P}}{x}$, and b) a collection of
edges that are square equivalent to $\Delta$.} 
\label{hyperplane}
\end{figure}

We describe combinatorial hyperplanes in $\kbpw$ with help from an example.  First, fix a $1$-cell
in $\kbpw$, such as the one in Figure \ref{hyperplane}a), which we'll call $\Delta$.

We consider a small number of $1$-cells that are square equivalent to $\Delta$ (there are infinitely many such
$1$-cells for this $\Delta$).  These are the vertical edges in Figure \ref{hyperplane}b).  We denote these edges
$\Delta_1$, $\Delta_2 = \Delta$, $\Delta_3$, and $\Delta_4$, reading from left to right.  Note the interpretation
of simple square equivalence in terms of diagrams:  for $i=1,2,3$, $\Delta_{i} \sim \Delta_{i+1}$ since $\Delta_{i}$
can be obtained from $\Delta_{i+1}$ by removing a maximal shaded transistor and all of its bottom wires
from $\Delta_{i+1}$, or the reverse, i.e., $\Delta_{i+1}$ can be obtained in the same way from $\Delta_{i}$.
This observation is general, and holds true in all of the complexes $\kpw$, $\kapw$, and $\kbpw$, for all $\mathcal{P}$
and $w$, and indeed follows easily from the definition of the $2$-cells in a diagram complex.  We record this 
in a lemma.

\begin{lemma}
Let $\Delta'$, $\Delta''$ be $1$-cells in $\kbpw$.  The following statements are equivalent:
\begin{enumerate}
\item $\Delta' $ and $\Delta''$ are simple square equivalent; 
\item There is some maximal shaded transistor $T$ in $\Delta'$ such that $\Delta''$ is
the result of removing $T$ and all of its bottom wires from $\Delta'$ (or the reverse statement
is true, with $\Delta'$ and $\Delta''$ reversing roles).
\qed
\end{enumerate}
\end{lemma}

Fix a $1$-cell $\Delta \subseteq \kbpw$.  Let $H_{\Delta}$ denote the combinatorial hyperplane corresponding to
$\Delta$.  Let $\mathcal{T}_{\Delta}$ denote the collection of transistors  of $\Delta$.
 Let $T$ denote the (unique) white transistor in $\mathcal{T}_{\Delta}$.  Consider the collection
$M_{\Delta} = \{   T' \mid   T' \mathrm{~is~ a~ transistor~in~} \Delta; \, \, T' \leq T \}$
(the inequality sign refers to the partial order on transistors).  We can associate to this collection of transistors
a vertex $min \left( H_{\Delta} \right)$, called the \emph{minimal vertex of $H_{\Delta}$}.  Simply remove all transistors
in $\mathcal{T}_{\Delta} - M_{\Delta}$ along with their bottom wires, and then shade the white transistor.  The result is
necessarily a vertex by Lemma \ref{bccc}.  It is clear that $min \left( H_{\Delta} \right)$ depends only on the hyperplane
$H_{\Delta}$.

For example, if $\Delta = \Delta_1 , \Delta_2 , \Delta_3,  ~ \mathrm{or} ~ \Delta_4$ from Figure \ref{hyperplane}b), then
$min \left( H_{\Delta} \right)$ is the vertex at the top of $\Delta_2$.  We note one property of minimal vertices:  a vertex
is minimal if and only if it contains a unique maximal transistor.  If a vertex $\Delta$ 
has a unique maximal transistor, 
then the hyperplane $H_{\Delta}$ corresponding to $\Delta$ is the square equivalence class of the edge obtained by painting the
maximal transistor white. 

We are interested in $min \left( H_{\Delta} \right)$ because of the following lemma.  In all that follows, 
we let our basepoint
$\ast$ be the unique vertex in $\kbpw$ having no transistors. 

\begin{lemma} \label{sep}
Let $H$ be a hyperplane in $\kbpw$.  If $\overline{\Delta}$ is an arbitrary vertex
in $\kbpw$, then $\overline{\Delta}$ and the basepoint $\ast$ lie in different components of
$\kbpw - H$ if and only if $min(H) \leq \overline{\Delta}$.
\end{lemma}
 
\begin{proof}
Let $\overline{\Delta}$ be a vertex in $\kbpw$; let $\mathcal{T}_{\overline{\Delta}}$ denote the collection 
of transistors in $\overline{\Delta}$.  Choose a function $\alpha : \mathcal{T}_{\overline{\Delta}}   
\rightarrow \left\{ 1, \ldots , \left| \mathcal{T}_{\overline{\Delta}} \right| \right\}$ satisfying:
\begin{enumerate}
\item $\alpha$ is one-to-one;
\item if $T_1 < T_2$, then $\alpha \left( T_1 \right) < \alpha \left( T_2 \right)$.
\end{enumerate}
We associate a sequence of vertices
$\ast = \overline{\Delta}_{0}, \overline{\Delta}_{1}, \ldots, 
\overline{\Delta}_{\mid \mathcal{T}_{\overline{\Delta}} \mid } = \overline{\Delta}$, where
$\overline{\Delta}_{i}$ is the (unique) vertex determined by $\alpha^{-1} \left( \{ 1, 2, \ldots, i \} \right)$
under the correspondence in Lemma \ref{bccc}.  It is not difficult to see that $\overline{\Delta}_{i}$ is 
connected to $\overline{\Delta}_{i+1}$
by a unique edge for $i=0, 1, \ldots, \left| \mathcal{T}\left( \overline{\Delta} \right) \right| - 1 $.  We let
$p_{\alpha}$ denote the edge-path consisting of these edges.

\begin{figure}[!h] 
\begin{center}
\includegraphics{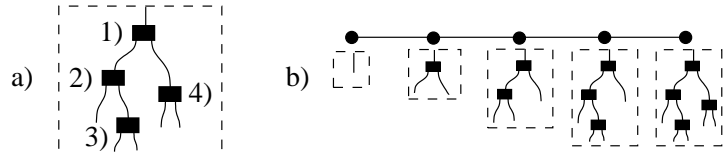}
\end{center}
\caption{a) A labelling $\alpha$ of the 
picture $\overline{\Delta}$, and b) the associated edge-path in $\kt{b}{\mathcal{P}}{x}$.} 
\label{edgepath}
\end{figure}

For example, Figure \ref{edgepath}a) shows a picture with a numbering $\alpha$ of its transistors, along with the corresponding 
edge-path ( Figure \ref{edgepath}b) ).

We claim that $p_{\alpha}$ is a geodesic in the $1$-skeleton $\kbpw^{1}$.
Suppose that $p$ is an arbitrary edge-path connecting $\ast$ to $\overline{\Delta}$; let
$\ast = \overline{\Delta}_{0}', \overline{\Delta}_{1}', \ldots, \overline{\Delta}_{m}' = \overline{\Delta}$
be the vertices lying along the path $p$, listed in the order they are visited.  It is clear from the definition
of edges in $\kbpw$ that $\overline{\Delta}_{i+1}'$ is obtained from $\overline{\Delta}_{i}'$ $( 0 \leq i \leq m-1 )$
by either removing a maximal transistor from the bottom of $\overline{\Delta}_{i}'$, or adding a new maximal transistor
to $\overline{\Delta}_{i}'$.  It immediately follows from this that 
$\ell(p) \geq \left| \mathcal{T}_{\overline{\Delta}} \right|$.  This proves the claim.

Suppose that $min \left( H_{\Delta} \right) \leq \overline{\Delta}$.  Lemma \ref{bccc} implies that
$min \left( H_{\Delta} \right)$ corresponds to an initial collection $\mathcal{T}$ of transistors.
Suppose that $\left| \mathcal{T} \right| = n$.  It follows that we can define 
$\alpha : \mathcal{T}_{\overline{\Delta}} \rightarrow \left\{ 1, \ldots, n, \ldots, 
\left| \mathcal{T}_{\overline{\Delta}} \right|  \right\}$ in such a way that
$\alpha_{\mid_{\mathcal{T}}} : \mathcal{T} \rightarrow \{ 1, \ldots, n \}$ is another labelling function satisfying
(1) and (2) above.  In this case, $\overline{\Delta}_{n} = min\left( H_{\Delta} \right)$ and 
$\overline{\Delta}_{n-1}$ is the vertex obtained by removing the (unique) maximal transistor in
$\overline{\Delta}_{n}$.  It immediately follows that the edge $\left[ \overline{\Delta}_{n-1} , \overline{\Delta}_{n} \right]$
is a member of the combinatorial hyperplane $H_{\Delta}$.  Since a geodesic edge-path $p_{\alpha}$ from $\ast$ to $\overline{\Delta}$
crosses $H_{\Delta}$, Theorem \ref{Sageev}(3) implies 
that $\ast$ and $\overline{\Delta}$ lie on opposite sides of $H_{\Delta}$, proving one direction.

Conversely, suppose that $\overline{\Delta}$ and $\ast$ are separated by the hyperplane $H_{\Delta}$.  
Theorem \ref{Sageev}(3) says 
that a geodesic edge-path $p$ crosses precisely the hyperplanes separating
the initial vertex of $p$ from the terminal vertex of $p$.  It follows that, for some
$k \in \left\{ 0, \ldots, \left| \mathcal{T}_{\overline{\Delta}} \right| - 1 \right\}$, the edge
$\left[ \overline{\Delta}_{k} , \overline{\Delta}_{k+1} \right]$ represents the hyperplane $H_{\Delta}$.
Under the correspondence in Lemma \ref{bccc}, $\overline{\Delta}_{k}$ corresponds to a collection
$\mathcal{T}_{\overline{\Delta}_{k+1}} - \{ T \}$ of transistors in $\overline{\Delta}_{k+1}$, 
where $T \in \mathcal{T}_{ \overline{\Delta}_{k+1}}$.  The edge
$\left[ \overline{\Delta}_{k} , \overline{\Delta}_{k+1} \right]$ can be described in terms of pictures as follows:
draw $\overline{\Delta}_{k+1}$, but leave the transistor $T$ unshaded.  According to the definition, we obtain
$min \left( H_{\Delta} \right)$ by shading $T$, and then taking the picture corresponding to
$\left\{ T' \in \mathcal{T}_{\overline{\Delta}_{k+1}} \mid T' \leq T \right\}$.  It is thus clear that
$min \left( H_{\Delta} \right) \leq \overline{\Delta}_{k+1} \leq \overline{\Delta}$.
\end{proof}  

\begin{proposition} \label{ord}
Let $H_{1}^{+}$ and $H_{2}^{+}$ be two positive half-spaces in $\kbpw$.  Let $\Delta_{1}$ and
$\Delta_{2}$ be their minimal vertices.
\begin{enumerate}
\item $H_{2}^{+} \leq H_{1}^{+}$ if and only if $\Delta_{2} \leq \Delta_{1}$.
\item $H_{1}^{+} \cap H_{2}^{+} \neq \emptyset$ if and only if $\{ \Delta_{1} , \Delta_{2} \}$
has an upper bound in $\kbpw^{0}$.
\end{enumerate}
\end{proposition}

\begin{proof}
(1) $( \Rightarrow )$ Suppose $H_{2}^{+} \leq H_{1}^{+}$.  This means that $H_{1}^{+} \subseteq H_{2}^{+}$.  Thus,
every vertex $\Delta \in H_{1}^{+}$ is separated from $\ast$ by $H_{2}$.  In particular, $\Delta_{1}$ is so separated
from $\ast$.  By the previous lemma, $\Delta_{2} \leq \Delta_{1}$.

$( \Leftarrow )$ Suppose $\Delta_{2} \leq \Delta_{1}$.  It is sufficient to check the inclusion
$H_{1}^{+} \subseteq H_{2}^{+}$ on vertices.  If $\Delta$ is a vertex in $H_{1}^{+}$, then
$\Delta_{1} \leq \Delta$.  It follows that $\Delta_{2} \leq \Delta$, so $\Delta \in H_{2}^{+}$ by the previous lemma.

(2) Both directions are immediate consequences of the previous lemma.
\end{proof}

We now obtain the desired characterization of profiles in terms of pictures.

\begin{theorem} \label{profile}
Let the basepoint $\ast \in \kt{b}{\mathcal{P}}{w}$ be the unique vertex having no transistors. 
Let $\Delta$ be an infinite $(w,\ast)$-picture over the semigroup presentation $\mathcal{P}$, 
i.e.,  a picture in the sense of Section \ref{diagram}, except that the transistor and wire sets
are countably infinite.  Let us suppose as well that $\Delta$ satisfies the following conditions:
\begin{enumerate}
\item For any transistor 
$T \in \mathcal{T}_{\Delta}$, the set $( -\infty, T] = \{ T' \in \mathcal{T}_{\Delta} \mid
T' \leq T \}$ is finite;
\item There are no maximal elements in the set $\mathcal{T}_{\Delta}$ of transistors in $\Delta$; 
\item No wire is attached to the bottom of the frame of $\Delta$.
\end{enumerate}
The picture $\Delta$ determines a unique profile, i.e., a 
non-empty collection of positive half-spaces $\mathcal{H}_{\Delta}$ in $\kbpw$  
satisfying       
properties (1)-(3) in Proposition \ref{fund}.  Conversely, a profile in $\kbpw$ determines a unique
infinite picture $\Delta$ satisfying properties (1)-(3) above.

The indicated correspondences are mutually inverse.
\end{theorem}

\begin{proof}
Suppose that $\Delta$ is an infinite $(w, \ast)$-picture over the semigroup presentation $\mathcal{P}$ satisfying
the properties above.  The transistors of $\Delta$ are in one-to-one correspondence with a collection of hyperplanes in
the following way.  Let $T$ be a transistor of $\Delta$; we consider the collection $(- \infty, T]$ of all transistors in
$\Delta$ which are less than or equal to $T$ in the partial order on transistors.  By Lemma \ref{bccc} and 
our assumption that $(- \infty, T]$ is finite, 
this collection of transistors
corresponds to a unique vertex, and this vertex is the minimal vertex of a unique hyperplane $H_{T}$.  Note that Lemma \ref{bccc}
also implies that the correspondence between transistors and hyperplanes is one-to-one.

We consider the properties of the collection $\mathcal{H}_{\Delta} = \left\{ H_{T} \mid T \in \mathcal{T}_{\Delta} \right\}$.  First, let
$H_{T_{1}}$, $\ldots$, $H_{T_{n}}$ be hyperplanes in $\mathcal{H}_{\Delta}$.  Consider the collection 
$\mathcal{T}' = \left\{ T \in \mathcal{T}_{\Delta} \mid T \leq T_{i} \mathrm{~for~some~} i \in \{ 1, \ldots, n \} \right\}$.  By 
Lemma \ref{bccc}, the collection $\mathcal{T}'$ corresponds to a vertex $\Delta_{\mathcal{T}'}$.  Moreover, we have
that $( -\infty, T_{i}] \subseteq \mathcal{T}'$, for $i=1, \ldots, n$, from which it follows that 
$min \left( H_{T_{i}} \right) \leq \Delta_{\mathcal{T}'}$, for $i=1, \ldots, n$.  This, in turn, implies that
$\Delta_{\mathcal{T}'} \in H_{T_{1}}^{+} \cap \ldots \cap H_{T_{n}}^{+}$, by Lemma \ref{sep}.  Thus property (1)
from Proposition \ref{fund} holds.

If $H_{T} \in \mathcal{H}_{\Delta}$, then, by the assumption that $T$ is not maximal, there is some 
$T_{1} \in \mathcal{T}_{\Delta}$ such that
$T < T_{1}$.  It follows from this that $(-\infty, T] \subsetneq ( -\infty, T_{1}]$; this implies that the
vertices $min \left( H_{T} \right)$, $min \left( H_{T_1} \right)$ under the correspondence from Lemma \ref{bccc}
satisfy $min \left( H_{T} \right) < min \left( H_{T_1} \right)$.  By Lemma \ref{ord}, $H_{T}^{+} < H_{T_1}^{+}$, 
so property (2) from Proposition \ref{fund} holds. 
   
Checking Property (3) from Proposition \ref{fund} is an easy exercise using the properties of the correspondence in Lemma \ref{bccc}.

Conversely, suppose that $\mathcal{H}$ is a non-empty collection of hyperplanes in $\kbpw$ satisfying properties (1)-(3) of Proposition \ref{fund}.
Choose a finite collection of
hyperplanes $H_{T_{1}}, H_{T_{2}}, \ldots, H_{T_{n}}$.
By property (1) from 
Proposition \ref{fund}, $H_{T_{1}}^{+} \cap H_{T_{2}}^{+} \cap \ldots \cap H_{T_{n}}^{+} \neq \emptyset$.  
This implies that there is some vertex $\hat{\Delta}$ in the latter intersection, which means, by Lemma \ref{sep}, that 
$min \left( H_{T_{i}} \right) \leq \hat{\Delta}$ for $i = 1, \ldots, n$.  Since the collection
$\left\{ min \left( H_{T_{i}} \right) \mid i = 1, \ldots, n \right\}$ has an upper bound, Lemma 3.2(2) of \cite{Me2}
implies that it has a least upper bound.  
Thus, we've shown that any finite collection of minimal vertices for hyperplanes in
$\mathcal{H}$ has a least upper bound.  (This least upper bound is a ``union" of the labels for these 
vertices, in an appropriate sense.
Note that it won't in general be a minimal vertex itself.) 

Let $\Delta_{1}, \Delta_{2}, \ldots, \Delta_{n}, \ldots $ be the sequence consisting of all minimal vertices for hyperplanes
in $\mathcal{H}$.  Since any finite collection of these hyperplanes has a least upper bound, we can identify the direct  
limit of this sequence with an infinite diagram $\Delta$.  It is clear that $\Delta$ has properties (1) and (3) from the statement of the Theorem. 
Property (2) follows easily from the fact that the collection $\mathcal{H}$ satisfies (2) from Proposition \ref{fund}.

We leave the final statement as an exercise.
\end{proof}      

We will sometimes require a lemma which gives a necessary condition
on the open cube $C$ through which a geodesic ray $c$ with profile $P(c)$
can travel.  The condition involves the largest vertex $\Delta_{C}$ in the
closure $\overline{C}$, which always exists, and can be obtained by
shading each transistor in the picture representative for $C$
(see for instance Figure \ref{cubepic}, from Section \ref{diagram}).

\begin{lemma} \label{trav}
Let $c: [0, \infty) \rightarrow \kbpw$ be a geodesic ray issuing from the base vertex
$\ast$.  Let $\Delta$ be the infinite picture representing the profile of $c$.  Let $C$
be an open cube of $\kbpw$ satisfying $ \left( \mathrm{Im} \, c \right) \cap C \neq \emptyset$.

If $\Delta_{C}$ is the largest vertex in $\overline{C}$, then $\Delta_{C} \leq \Delta$.
\end{lemma}

\begin{proof}
Let $\Delta_{C}$ be the largest vertex of $C$.  Consider the collection of all maximal transistors
$T_{1}, \ldots, T_{n}$ in $\Delta_{C}$; let $\Delta_{T_{i}}$ be the unique vertex determined by
$( -\infty, T_{i}]$ under the correspondence from Lemma \ref{bccc}.  Note that each $\Delta_{T_{i}}$
is the minimal vertex  of a hyperplane $H_{T_{i}}$, and all of the hyperplanes
$H_{T_{1}}, \ldots, H_{T_{n}}$ are distinct.

Let $\widehat{\Delta}_{C}$ denote the representative for $C$, which is a picture consisting of shaded and
unshaded transistors, as in Figure \ref{cubepic}.  Some of the transistors $T_i$ are shaded in $\widehat{\Delta}_{C}$;
others are unshaded in $\widehat{\Delta}_{C}$.  

If $T_i$ is shaded, then $\Delta_{T_{i}} \leq \Delta_{C}$, so that $H_{T_{i}}$ separates $\Delta_{C}$ from $\ast$.
Moreover, $H_{T_{i}}$ doesn't pass through $C$, since the hyperplanes $H$ satisfying $H \cap C \neq \emptyset$ are the precisely
the collection of all $H_{T_{j}}$ such that $T_j$ is unshaded in $\widehat{\Delta}_{C}$.  It follows that each point $x$ in 
$\left( \mathrm{Im} \, c \right) \cap C$ can be connected to $\Delta_{C}$ without crossing $H_{T_{i}}$; therefore $x \in H^{+}_{T_{i}}$,
so $H^{+}_{T_{i}} \in P(c)$.

If $T_i$ is unshaded in $\widehat{\Delta}_{C}$, then $H_{T_{i}} \cap C \neq \emptyset$.  It follows from Lemma \ref{crossing}(2) that
$H^{+}_{T_{i}} \in P(c)$.

The correspondence of Theorem \ref{profile} implies that $\Delta_{T_{1}}, \ldots, \Delta_{T_{n}} \leq \Delta$.  This means that the least
upper bound $\widetilde{\Delta}$ of $\left\{ \Delta_{T_{1}}, \ldots, \Delta_{T_{n}} \right\}$ exists, and satisfies $\widetilde{\Delta} \leq \Delta$.
But clearly $\widetilde{\Delta} = \Delta_{C}$, since the latter vertex is an upper bound of $\left\{ \Delta_{T_{1}}, \ldots, \Delta_{T_{n}} \right\}$,
and any proper initial subset subset of $\mathcal{T}_{\Delta_{C}}$ would fail to contain at least one of the transistors $T_{1}, \ldots, T_{n}$.
\end{proof}

\subsection{The Action on Profiles}

\begin{proposition} \label{action}
There is a well-defined action of $\dbpw$ on the set of all profiles.  If $\Delta$ is
a profile and $\Delta_{1} \in \dbpw$, then $\Delta_{1} \ast \Delta$ can be computed as
follows.  First, form the concatenation $\Delta_{1} \circ \Delta$ and remove all dipoles.
Second, remove all maximal transistors from the resulting infinite diagrem, until no maximal
transistors remain.  The result is $\Delta_{1} \cdot \Delta$.
\end{proposition}
  
\begin{proof}
Let $\mathcal{H}$, $\mathcal{H}'$ be two collections of positive half-spaces in 
$\kbpw$.  We write $\mathcal{H} \sim \mathcal{H}'$ and $[\mathcal{H}] = [\mathcal{H}']$
if $\mathcal{H}$ and $\mathcal{H}'$ are \emph{cofinal}, i.e., if for any
$H_{1}^{+} \in \mathcal{H}$, there exists
$\left( H'_{1} \right)^{+} \in \mathcal{H}'$ such that $H_{1}^{+} \leq \left( H_{1}' \right)^{+}$ and
for any $\left( H_{2}' \right)^{+} \in \mathcal{H}'$, there exists $H_{2}^{+} \in \mathcal{H}$
such that $\left( H_{2}' \right)^{+} \leq H_{2}^{+}$.  It is fairly clear that $\sim$ is an
equivalence relation on the set of all collections of positive half-spaces in $\kbpw$.

The group $\dbpw$ doesn't act in an obvious way on the set of equivalence classes, since a group
element $\overline{\Delta} \in \dbpw$ doesn't necessarily map a positive half-space to a positive
half-space.  Indeed, $\overline{\Delta} \cdot H^{+}$ is a negative half-space if and only if
$\ast \in \overline{\Delta} \cdot H^{+}$, that is, if and only if $\overline{\Delta}^{-1} \cdot \ast \in H^{+}$.
Now $\overline{\Delta}^{-1} \cdot \ast \in H^{+}$ if and only if $min \left( H^{+} \right) \leq \overline{\Delta}^{-1} \cdot \ast$.
There are only finitely many vertices $\Delta'$ satisfying $\Delta' \leq \overline{\Delta}^{-1} \cdot \ast$ (all determined by
the correspondence in Lemma \ref{bccc}).  It follows from this that a given $\overline{\Delta} \in \dbpw$ maps at most finitely 
many positive half-spaces to negative ones.

We obtain an action on profiles in the following way.  
Identify a profile $\Delta$ with the unique equivalence class $\left[ \mathcal{H} \right]$
such that $\Delta \in \left[ \mathcal{H} \right]$.  For a given $\overline{\Delta} \in \dbpw$, choose a collection
$\mathcal{H}' \in \left[ \mathcal{H} \right]$ such that $\overline{\Delta} \cdot H^{+}$ is a positive half-space,
for each $H^{+} \in \mathcal{H'}$.  
It is possible to do this since each $\mathcal{H}' \in \left[ \mathcal{H} \right]$ 
is necessarily infinite.  We define
$\overline{\Delta} \ast \Delta$ to be $\left[ \overline{\Delta} \cdot \mathcal{H}' \right]$.  It is not difficult to see
that $\left[ \overline{\Delta} \cdot \mathcal{H}' \right]$ contains a 
unique profile, and that the definition of $\ast$ doesn't depend
upon the choice of $\mathcal{H}' \in \left[ \mathcal{H} \right]$.  It follows that $\ast$ is an action on profiles.

We claim two things:  first, that $\Delta \ast P(c) = P ( \Delta \ast c )$, for any geodesic ray
$c: [0, \infty) \rightarrow \kbpw$ issuing from $\ast$; second, that the action $\ast$ from the previous paragraph
has the description promised in the statement of the proposition.

Recall the definition of the action $\ast: \dbpw \times \partial \kbpw \rightarrow \partial \kbpw$ on the
space at infinity.  If $c \in \partial \dbpw$, then, for any $\overline{\Delta} \in \dbpw$, 
$\overline{\Delta} \cdot c$ is simply the translate of $c$ by the usual action of $\overline{\Delta}$ on
$\kbpw$.  The ray $\overline{\Delta} \ast c$ is the unique ray issuing from $\ast$ and asymptotic to $\Delta \cdot c$,
the existence of which is guaranteed by Proposition \ref{parallel}.

We now prove the first claim.  Let $c \in \partial \kbpw$ and let $\overline{\Delta} \in \dbpw$.  We choose
some cofinal subset $\mathcal{H} \subseteq P(c)$ such that $\overline{\Delta} \cdot H^{+}$ is a positive half-space,
for any $H^{+} \in \mathcal{H}$.  It is clear that $\overline{\Delta} \cdot c$ crosses each of the hyperplanes
in $\overline{\Delta} \cdot \mathcal{H}$, since $c$ crosses each of the hyperplanes in $\mathcal{H}$.  Since
$\overline{\Delta}$ maps positive half-spaces in $\mathcal{H}$ to positive half-spaces, $\overline{\Delta} \cdot c$
intersects each $\overline{\Delta} \cdot H^{+} \in  \overline{\Delta} \cdot \mathcal{H}$ in an open ray, just as $c$
intersects each $H^{+} \in \mathcal{H}$ in an open ray.  By Lemma \ref{monotone},
$d_{\overline{\Delta} \cdot H} \left( \overline{\Delta} \cdot c(t) \right) \rightarrow \infty$ as 
$t \rightarrow \infty$ for any $H$ such that $H^{+} \in \mathcal{H}$.  It follows from this, first, that
$d_{\overline{\Delta} \cdot H} \left( \overline{\Delta} \ast c(t) \right)$ also goes to infinity as $t \rightarrow \infty$,
and second, that $\left( \overline{\Delta} \ast c \right) (t) \in \overline{\Delta} \cdot H^{+}$ for $t$ sufficiently large.
This implies that $\overline{\Delta} \cdot \mathcal{H} \subseteq P \left( \overline{\Delta} \ast c \right)$.

Next, we need to show that, for any positive half-space $H^{+} \in P \left( \overline{\Delta} \ast c \right)$, there
is $\overline{H}^{+} \in \overline{\Delta} \cdot \mathcal{H}$ such that $H^{+} \leq \overline{H}^{+}$.  Choose a sequence
of positive half-spaces $H^{+} = H^{+}_{0} < H^{+}_{1} < H^{+}_{2} < \ldots$ in $P \left( \overline{\Delta} \ast c \right)$.
By Lemma \ref{monotone}, we know that $d_{H_{i}} \left( ( \overline{\Delta} \ast c ) (t) \right) \rightarrow \infty$ as
$t \rightarrow \infty$, for any $i \in \{ 0, 1, \ldots, n, \ldots \}$.  By the definition of $P \left( \overline{\Delta} \ast c \right)$,
for any $i$, $\left( \overline{\Delta} \ast c \right) (t) \in H^{+}_{i}$ for $t$ sufficiently large.  It now follows from the fact that
$\overline{\Delta} \cdot c$ and $\overline{\Delta} \ast c$ are asymptotic that, for any $i$, $\left( \overline{\Delta} \cdot c \right)(t) \in H^{+}_{i}$
for $t$ sufficiently large.  This implies that $\overline{\Delta} \cdot c$ crosses at least one of the $H_i$, and therefore all $H_j$ for $j$ sufficiently large,
for otherwise
$$\left( \overline{\Delta} \cdot c \right) \left( [0, \infty) \right) \subseteq \bigcap_{i=0}^{\infty} H^{+}_{i} = \emptyset.$$ 
We choose some $H_{j}^{+}$ large enough that $\overline{\Delta}^{-1} \cdot H_{j}^{+}$ is a positive half-space.
Since $c$ crosses $\overline{\Delta}^{-1} \cdot H_{j}$, it follows that $\overline{\Delta}^{-1} \cdot H_{j}^{+} \in P(c)$.
Since $\mathcal{H}$ is cofinal in $P(c)$, there is $\hat{H}^{+} \in \mathcal{H}$ such that 
$\overline{\Delta}^{-1} \cdot H_{j}^{+}
\leq \hat{H}^{+}$.  This implies that $H_{j}^{+} \leq \overline{\Delta} \cdot \hat{H}^{+} \in \overline{\Delta} \cdot 
\mathcal{H}$.  Now
we've shown that $H^{+} \leq H_{j}^{+} \in \overline{\Delta} \cdot \mathcal{H}$.

It follows that $P \left( \overline{\Delta} \ast \mathcal{H} \right) = \left[ \overline{\Delta} \cdot \mathcal{H} \right]$ under the identification
of $\left[ \overline{\Delta} \cdot \mathcal{H} \right]$ with a profile.  We've thus shown that $P \left( \overline{\Delta} \ast c \right) =
\overline{\Delta} \ast P(c)$.  It immediately follows from this that the action of $\dbpw$ on profiles is well-defined:
if $c, c'$ have the same profile, then so also do $\overline{\Delta} \ast c$ and $\overline{\Delta} \ast c'$.

Now we prove the second claim.  Let $\Delta$ be an infinite picture representing a profile and let
$\overline{\Delta} \in \dbpw$.  We choose a cofinal collection of transistors $\mathcal{T}'$ in $\Delta$ (which are identified with positive half-spaces
by the correspondence in Theorem \ref{profile}) such that no transistor in $\mathcal{T}'$ forms a dipole in the
concatenation $\overline{\Delta} \circ \Delta$.  The above description of the action implies that
$\overline{\Delta} \ast \Delta$ is the collection of all positive half-spaces such that
$H^{+} \leq \overline{\Delta} \cdot H_{T'}^{+}$, for some $T' \in \mathcal{T}'$.  This collection of positive half-spaces
may be identified with the collection of transistors $T$ in the reduced concatenation $\overline{\Delta} \circ \Delta$
satisfying $T \leq T'$, for some $T' \in \mathcal{T}'$.  By the cofinality of $\mathcal{T}'$ in $\Delta$, these transistors $T$ are precisely those
for which there exists an infinite sequence $T = T_0 < T_1 < \ldots < T_n < \ldots$ where each $T_i$ is a transistor of
$\overline{\Delta} \circ \Delta$.  The second claim follows.
\end{proof}


\section{Fixed Profiles under the Actions of $F$, $T$, and $V$} \label{main}


The results of the previous section largely
reduce the problem of finding 
fixed points in $\partial \kt{}{\mathcal{P}}{x}$,  $\partial \kt{a}{\mathcal{P}}{x}$, 
and $\partial \kt{b}{\mathcal{P}}{x}$ ( where $\mathcal{P} = \langle x \mid x = x^{2} \rangle$) to   
the algebraic problem of finding globally fixed profiles.  The latter problem is quite
easy; we give a complete classification of fixed profiles for $F$, $T$, and $V$ in
this section.

\subsection{Conventions}

We fix some conventions for portraying profiles (and pictures) over the semigroup
presentation $\langle x \mid x = x^{2} \rangle$.  

\begin{figure} [!h] 
\begin{center}
\includegraphics{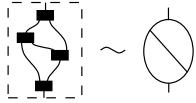}
\end{center}
\caption{A convention for drawing pictures over the semigroup presentation
$\mathcal{P} = \langle x \mid x=x^{2} \rangle$.  On the left, we have a picture drawn in 
the usual fashion; on the right is its equivalent. } 
\label{gen}
\end{figure}

First, we draw every transistor in a
picture or profile as a point and omit the frame, so that, for instance, the element $x_0 \in F$
(depicted as an ordinary picture on the left) looks like the right half of Figure \ref{gen}.

We also need some conventions that will allow us to portray infinite pictures using a finite amount
of space.
We draw an empty dot at the end of a wire to indicate that the bottom of the wire doesn't
connect to any transistor.  A solid dot at the end of a wire indicates that the wire
connects to the top of some transistor.  This transistor may be either an $(x,x^{2})$-transistor
or an $(x^{2}, x)$-transistor.  
A wire with no dot at the end may connect to a transistor or not; we make no assumption one way or the
other.

Finally, we let $T_{m}$ denote the full ordered rooted binary tree of depth $m$.  We let $\dot{T}_{m}$
denote the full ordered rooted binary tree of depth $m$, where each leaf ends in a solid dot.
Thus, $\dot{T}_m$ is the unique ordered rooted binary tree having $2^{m}$ dotted leaves, each at distance $m$ from the root.
The dot on each leaf (wire) implies that each connects to the top of some transistor.  The picture $T_{m}$ is the same tree, but without
the dots on the leaves.  
Thus, no particular wire in $T_m$ which corresponds to a leaf necessarily leads 
to the top of a transistor.  Notice, however,
that in a profile 
at least one of the leaves beneath a given degree $3$ vertex must attach to the top of a transistor, since there are no maximal transistors
in a profile (by Theorem \ref{profile}). 

\subsection{Thompson's Group $F$}

\begin{figure} [!h] 
\begin{center}
\includegraphics{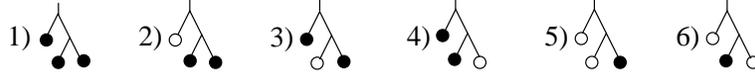}
\end{center}
\caption{The six cases.}
\label{cases}
\end{figure}

Let $\Delta$ be a profile of $\kpw$ which is fixed by all of $F$.  Without loss of generality,
we can assume that $\Delta$ has one of the forms in Figure \ref{cases}.
(The only other possible cases are $1' - 6'$, which are obtained by reflecting $1 - 6$ across a vertical axis.  Note
that the resulting cases are not mutually exclusive.)

\subsubsection{The Even-numbered Cases}

\begin{figure} [!h] 
\begin{center}
\includegraphics{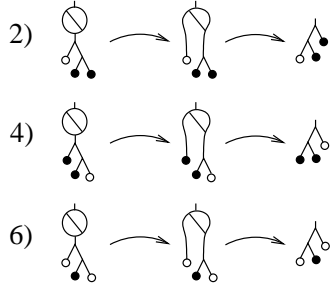}
\end{center}
\caption{The even-numbered cases.}
\label{evencases}
\end{figure}

We consider even-numbered cases first.  Let $x_0$ act on any profile $\Delta$ covered by Case $2$, $4$, or $6$.  The
results appear in Figure \ref{evencases}.
After cancelling dipoles twice, we arrive at the infinite pictures
in the column at the far right of the Figure.  We claim
that these infinite pictures 
necessarily contain no dipoles, no matter how the wires terminating in  black dots are connected
to transistors.  (It is clear also that these contain no maximal transistors.)

To prove the claim, first note that any dipole in the product $x_0 \cdot \Delta$ must be formed of one transistor
in $x_0$ and another in $\Delta$.  Thus, an infinite picture 
on the right side of Figure \ref{evencases} contains a dipole only if
one of the two pictured vertices of degree $3$ 
(both of which represent transistors from $x_0$) can form the top half of a dipole.

\begin{figure} [!h] 
\begin{center}
\includegraphics{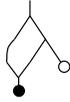}
\end{center}
\caption{An impossible profile.  No matter how we attach the black-dotted wire to the top of a
transistor, a dipole will be formed.} 
\label{cont}
\end{figure}

This is clearly impossible in Cases $2$ and $6$.  In Case $4$, it is enough to show
that the lower transistor cannot form the top half of a dipole. 
If we assume that it does, then the original profile $\Delta$ would have the form in Figure \ref{cont}.
No profile can have this form, however, since it is impossible to connect the black-dotted wire to a transistor
without forming a dipole, and $\Delta$ cannot contain dipoles.  This proves the claim. 

Finally, we compare the reduced profiles $x_0 \ast \Delta$ at the right in Figure \ref{evencases} with the originals in
Figure \ref{cases}. Since $\Delta$ is fixed by all of $F$, we must have that $x_0 \ast \Delta = \Delta$.  This is
impossible, as we easily see.  For instance, in Case $2$, the left wire dangling from the bottom of the topmost
transistor in $\Delta$ doesn't connect to a transistor, but the wire of the same description in $x_0 \ast \Delta$
does.  The other even cases are left as easy exercises.

\subsubsection{Cases 3 and 5}
\begin{figure} [!h] 
\begin{center}
\includegraphics{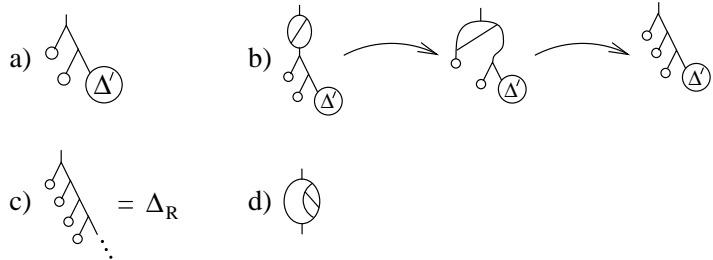}
\end{center}
\caption{a) The general form of a profile from Case 5.  b) The effect of letting $x_{0}^{-1}$ act.
c) The profile $\Delta_{R}$.  d)  The element $x_1 \in F$.} 
\label{sub5}
\end{figure}

If $\Delta$ is represented in Case $5$, then it must have the form in Figure \ref{sub5}a), where $\Delta'$ is another
profile.
If we let $x_{0}^{-1}$ act on $\Delta$, then after removing a dipole and an exposed transistor, we arrive at the
profile on the far right of Figure \ref{sub5}b).  A simple induction using the fact that $\Delta = x_{0}^{-1} \ast \Delta$
now shows $\Delta$ is the (unique) profile of the form depicted in Figure \ref{sub5}c).  It is easy to check that
$\Delta$ is fixed by all of $F$; indeed, it is enough to show that $\Delta$ is fixed by the generators $x_0$ and $x_1$.
We leave this verification as an exercise.

\begin{figure} [!h] 
\begin{center}
\includegraphics{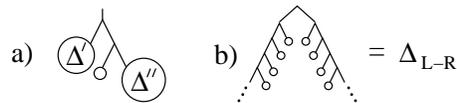}
\end{center}
\caption{a) The general form of a profile from Case 3.  b) The profile $\Delta_{L-R}$.}
\label{lr}
\end{figure}

If $\Delta$ is represented in Case $3$, then it must have the form in Figure \ref{lr}a), where $\Delta'$ and $\Delta''$
are profiles.
An argument similar to the one used for Case $5$ shows that $\Delta''$ has the form depicted in Figure \ref{sub5}c), and
$\Delta'$ is the result of reflecting $\Delta''$ across a vertical axis.  The details are left as an exercise.  It
follows  that $\Delta$ is the profile depicted in Figure \ref{lr}b), which is indeed fixed by both $x_0$ and $x_1$.

\begin{figure} [!h] 
\begin{center}
\includegraphics{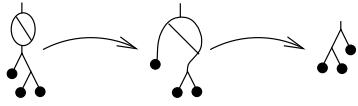}
\end{center}
\caption{The action of $x_0$ on a profile from Case 1.} 
\label{case1}
\end{figure}

\subsubsection{Case 1}

\begin{figure} [!h] 
\begin{center}
\includegraphics{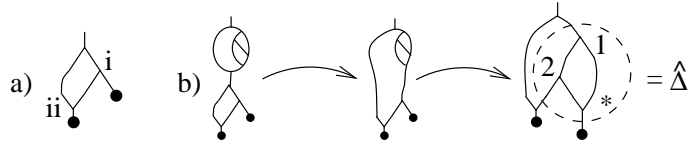}
\end{center}
\caption{a)  If we let $x_0$ act on this profile, the transistor labelled ii) will form half of a dipole.  b) The action
of $x_1$ on the profile from a).}
\label{subcases1}
\end{figure}

We now turn to Case $1$.  Let $x_0$ act on $\Delta$.
There are two subcases to consider:  either $\Delta$ has the form depicted in Figure \ref{subcases1}a) 
(and thus the infinite picture 
at the far right in Figure \ref{case1} contains a dipole) or the infinite picture 
at the far right in Figure \ref{case1} is reduced. 

 We
now rule out the first possibility using the fact that $\Delta$ is invariant under the action of $F$.  Let
$x_1$ act on $\Delta$; after reducing two dipoles we arrive at the infinite picture 
$\widehat{\Delta}$ on the far right of
Figure \ref{subcases1}b).  The transistors enclosed by the dotted circle 
were contributed by $x_1$, and any dipole in $\widehat{\Delta}$ would have to involve one of these three 
transistors.  Now note that, 
of these, only the transistor labelled $\ast$ could form half of a dipole; the others could not, even
after we cancel any dipole involving $\ast$.  Now we compare $\Delta$ and $\widehat{\Delta}$.  Under any isomorphism
between $\Delta$ and $\widehat{\Delta}$, the transistors labelled $1$ and $2$ in $\widehat{\Delta}$ must correspond
(respectively) to the transistors labelled i) and ii) in $\Delta$ (as depicted in Figure \ref{subcases1}a).  This is not possible,
since ii) is a $(x^{2}, x)$-transistor and $2$ is a $(x, x^{2})$-transistor.

\begin{figure} [!h] 
\begin{center}
\includegraphics{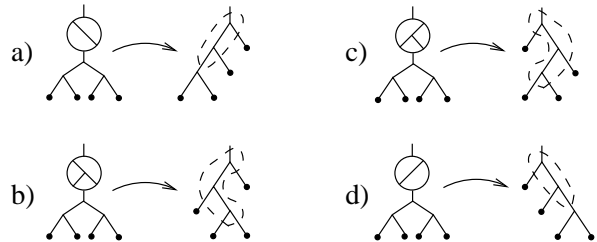}
\end{center}
\caption{The actions of a) $x_0$, b) $x_0 x_1 x_{0}^{-1}$, c) $x_1 x_{0}^{-1}$ and d) $x_{0}^{-1}$ on $\Delta$.}
\label{arg}
\end{figure}

It follows that we can assume that there are no dipoles in the 
infinite picture at the far right of Figure \ref{case1}.  After comparing
this profile with the profile $\Delta$, we can conclude that $\Delta$ has the form $\dot{T}_{2}$. 
We now multiply $\Delta$ by $x_0$, $x_0 x_1 x_{0}^{-1}$, $x_1 x_{0}^{-1}$ and $x_{0}^{-1}$.  The results are listed in
Figure \ref{arg}(a-d) (in the same order).
Dotted circles enclose the transistors that were contributed by the acting element.  If we knew that there were
no dipoles in the infinite 
pictures at the right, we could use the fact that all of them are equal to $\Delta$ in order to
conclude that $\Delta$ has the form  $\dot{T}_{3}$.  
\begin{figure} [!h]
\begin{center}
\includegraphics{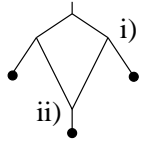}
\end{center}
\caption{This profile would form dipoles in cases b) and c) from Figure \ref{arg}.}
\label{lucky}
\end{figure}
The infinite pictures 
on the far right of a) and d) in Figure \ref{arg} are necessarily reduced.  The profiles in b) and c) will
be reduced unless $\Delta$ has the form in Figure \ref{lucky}b).

We now rule out the latter possibility.  Let $x_1$ act on $\Delta$.
\begin{figure} [!h]
\begin{center}
\includegraphics{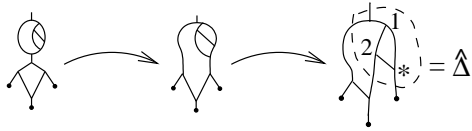}
\end{center}
\caption{The action of $x_1$ on the profile from Figure \ref{lucky}.}  
\label{lucky2}
\end{figure}
Any isomorphism between $\widehat{\Delta}$ and $\Delta$ must match the transistors labelled 1) and 2) with
the transistors labelled i) and ii), respectively.  This is impossible, since 2) is an $(x, x^{2})$-transistor and
ii) is an $(x^{2}, x)$-transistor.

\begin{figure} [!h]
\begin{center}
\includegraphics{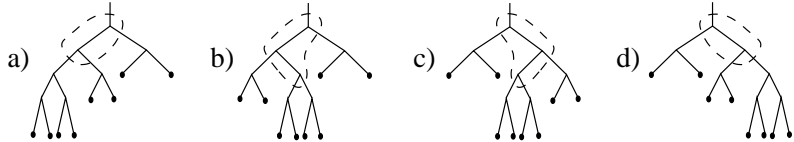}
\end{center}
\caption{The results of letting a) $x_0$, b) $x_0 x_1 x_{0}^{-1}$, c) $x_1 x_{0}^{-1}$, and d) $x_{0}^{-1}$ act on
$\dot{T}_{3}$.}
\label{deeper}
\end{figure}

It now follows that $\Delta$ has the form of $\dot{T}_{3}$.  We multiply $\Delta$ by $x_0$, $x_0 x_1 x_{0}^{-1}$, 
$x_1 x_{0}^{-1}$, $x_{0}^{-1}$.  The results appear in Figure \ref{deeper}.
Note that it is no longer possible for the transistors from the acting elements (circled) to form dipoles, so
all of the profiles in Figure \ref{deeper} are reduced.  Each of these profiles is equal to $\Delta$, since $\Delta$ is invariant
under the action of $F$.  It follows that $\Delta$ has the form of $\dot{T}_{4}$.

We now repeat this argument, letting the same four elements act on $\Delta$.  In this way, we conclude by induction that
$\Delta$ has the form of the full infinite binary $T_{\infty}$.  We write $\Delta = \Delta_{\infty}$.

We've proved the following theorem:
\begin{theorem}
Thompson's group $F$ fixes exactly four profiles:
$$ \Delta_{L}, \Delta_{R}, \Delta_{L-R}, ~and~ \Delta_{\infty}.$$
\qed
\end{theorem}
These profiles come from Cases $5'$, $5$, $3$, and $1$, respectively.

\subsection{Thompson's Groups $T$ and $V$}

This subsection is devoted to a proof of the following theorem:
\begin{theorem} \label{fixTV}
Thompson's groups $T$ and $V$ fix only the profile $\Delta_{\infty}$.
\end{theorem}
\begin{proof}
Let $\mathcal{P} = \langle x \mid x = x^{2} \rangle$.
Suppose that the group $\di{a}{\mathcal{P}}{x} \cong T$ fixes the profile $\Delta$.  It is not difficult
to see that $\Delta$ has the form in Figure \ref{TV}a), without loss of generality.
\begin{figure} [!h]
\begin{center}
\includegraphics{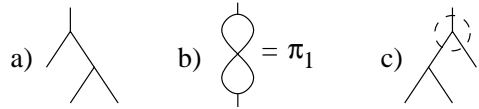}
\end{center}
\caption{a) The general form of a profile in the cubical complexes for $T$ and $V$; b) the acting
element $\pi_1$, which represents a half-rotation of the circle; c) the effect of letting $\pi_1$ act on
the profile from a).}       
\label{TV}
\end{figure}
We let $\pi_{1} \in \di{a}{\mathcal{P}}{x}$ act on $\Delta$; the result is portrayed in Figure \ref{TV}c), where
the circled transistor is contributed by $\pi_1$.  It follows that the picture
in c) is reduced.  Combining a) and c), which are both equivalent since $\pi_1 \ast \Delta = \Delta$, we have
that $\Delta$ has the form $T_{2}$.

Next, we claim that at least one of the four leaves of $T_{2}$ connects to the top of an
$(x, x^{2})$-transistor.  If not, then consider an $(x, x^{2})$-transistor $T$ in $\Delta$
which is distinct from the three $(x, x^{2})$-transistors in $T_{2}$, and minimal among $(x, x^{2})$-transistors
with this property.  The wire attached to the top of such a transistor could only lead up to
the bottom of an $(x^{2}, x)$-transistor $T'$.  This implies that $T$ and $T'$ form a dipole, which
contradicts the fact that $\Delta$ is reduced.

We therefore assume, without loss of generality, that $\Delta$ has the form in Figure \ref{wlog}a).
\begin{figure} [!h]
\begin{center}
\includegraphics{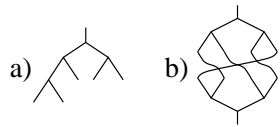}
\end{center}
\caption{a) Without loss of generality, $\Delta$ has this form.  b) The element $\pi_2$, which represents
a quarter-turn of the circle.}  
\label{wlog}
\end{figure}
Now we multiply $\Delta$ by $\pi_2$, $\pi_{2}^{2}$, and $\pi_{2}^{3}$, where $\pi_{2}$ is as in Figure \ref{wlog}b)    
to get the three profiles in Figure \ref{three}.
\begin{figure} [!h]
\begin{center}
\includegraphics{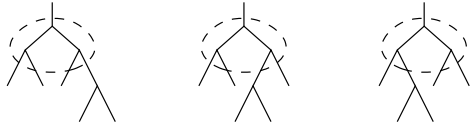}
\end{center}
\caption{The results of letting the powers of $\pi_2$ act on the profile from Figure \ref{wlog}a).} 
\label{three}
\end{figure}
Since the profiles from Figure \ref{three} and the profile from Figure \ref{wlog}a) are all equal to $\Delta$, it follows that $\Delta$ has the form
$T_3$.

We then argue, as before, that at least one of the eight leaves at the bottom of $T_3$ must be attached to the top of
an $(x, x^{2})$-transistor.  We then multiply $\Delta$ by $1$, $\pi_{3}$, $\pi_{3}^{2}$, $\pi_{3}^{3}$, $\ldots$,
$\pi_{3}^{7}$, where $\pi_{3}$ is a picture representing a one-eighth turn of the circle. 

If we compare the eight resulting profiles, we conclude that $\Delta$ is equivalent to $T_4$.  We can continue in a similar
way, and eventually conclude that $\Delta = \Delta_{\infty}$.

This proves the theorem in the case of $T$, and the proof for $V$ is the same word for word.
\end{proof}


\section{The Cases of $\Delta_{L}$, $\Delta_{R}$, and $\Delta_{L-R}$} \label{end}

\begin{figure} [!h]  
\begin{center}
\includegraphics{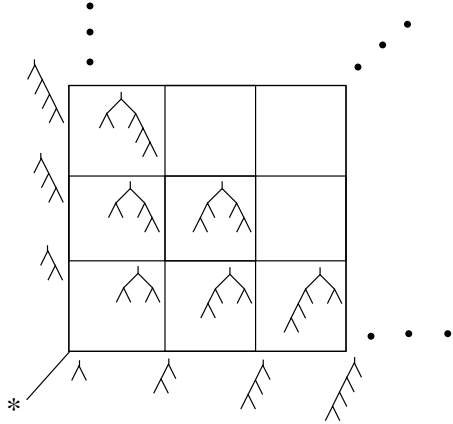}
\caption{A flat sector in $X_F$.  An integer lattice point $(m,n)$ ($m,n \geq 0$) corresponds to the 
tree $T_{m,n}$ having a root caret, $m$ carets dangling to the left, and $n$ carets dangling to the right.} 
\label{funky}
\end{center}
\end{figure}

Now we consider the geodesic rays $c \in \Delta_{L} \cup \Delta_{R} \cup \Delta_{L-R}$.  Consider 
first the profile $\Delta_{L-R}$.  By Lemma \ref{trav}, any geodesic ray $c$ in $\Delta_{L-R}$ is contained in
the subcomplex $K$ of $X_F \, \, ( = \dpw$, where $\mathcal{P} = \langle x \mid x = x^{2} \rangle$ and $w=x$). 
pictured in Figure \ref{funky}. 

Thus $K$ may be naturally identified with $\mathbb{R}^{2,+} \cup I$, where
$\mathbb{R}^{2,+} = \{ (x,y) \in \mathbb{R}^{2} \mid x,y \geq 0 \}$, $I$ is the unit interval, and 
$\mathbb{R}^{2,+} \cap I = \{ (0,0) \}$.

\begin{lemma} \label{Crisp} 
The inclusion $i: K \rightarrow X_F$ is an isometric embedding.
\end{lemma}

\begin{proof}
We appeal to Theorem 1(2) of \cite{CW}, which says:
If $X$ and $Y$ are finite dimensional CAT(0) cubical complexes and $\Phi: X \rightarrow Y$ is a cubical map, then  
the map $\Phi$ is an isometric embedding if and only if, for every vertex $v \in X$, the simplicial map between links
$\mathrm{Lk} (x,X) \rightarrow \mathrm{Lk}(\Phi(x), Y)$ induced by $\Phi$ is injective with image a full subcomplex
of $\mathrm{Lk}(\Phi(x),Y)$.  (We refer the reader to \cite{BH}, page 102 for a discussion of the link; Crisp and Wiest
define cubical maps on page 443 of \cite{CW}, and it is clear that the inclusion map is cubical.)

We consider the link of a vertex $T_{m,n}$, where $m,n > 1$, and leave the verifications for the other vertices as
an exercise.  
The link $\mathrm{Lk}(T_{m,n},K)$ is a square, i.e., the obvious one-dimensional simplicial complex consisting of $4$ vertices and $4$ edges.
This link will be embedded in $\mathrm{Lk}(T_{m,n},X_F)$ as a full subcomplex if and only if (1) there is no two-dimensional cube $C$ in
$X_F$ such that $T_{m,n-1}$ and $T_{m,n+1}$ are both vertices of $C$, and (2) there is no two-dimensional cube $C$ in $X_F$ such that
$T_{m-1, n}$ and $T_{m+1, n}$ are both vertices of $C$.        

We now check (1); the argument for (2) is similar.  
If there is such a cube $C$, then $C$ can be represented by a picture as in Figure \ref{cubepic}, consisting of $k$ shaded
and $2$ unshaded transistors.  The four corners of $C$ are labelled by pictures having $k$, $k+1$, $k+1$, and $k+2$ transistors.  It follows
that $k = m+n$, that $T_{m,n-1}$ is the result of leaving off both unshaded transistors, and that $T_{m,n+1}$ is the result of 
shading both transistors.  From this we get a contradiction, since the unshaded transistors of the cube $C$ must both be maximal, and there is
no way to remove two transistors that are both maximal in $T_{m,n+1}$ and arrive at $T_{m,n-1}$.
\end{proof}

\begin{theorem}
If $c$ is a geodesic ray in $X_F$ issuing from $\ast$ and $c \not \in \Delta_{\infty}$, then $c$ represents a point at infinity that is fixed by all
of $F$ if and only if $c \in \Delta_{L} \cup \Delta_{R} \cup \Delta_{L-R}$.  The subspace of $\partial X_F$
consisting of $\Delta_{L} \cup \Delta_{R} \cup \Delta_{L-R}$ is an arc of Tits length $\pi / 2$.
\end{theorem}

\begin{proof}
If $f: X \rightarrow Y$ is an isometric embedding between CAT(0) spaces, then the induced map $f_{\infty}: \partial X
\rightarrow \partial Y$ is an isometry, where the boundary is endowed with the angular metric (\cite{BH}, page 280).    
By the previous lemma, $K$ is isometrically embedded in $X_F$; by Lemma \ref{trav}, a geodesic ray $c$ issuing from $\ast$
represents $[c] \in \Delta_{L} \cup \Delta_{R} \cup \Delta_{L-R}$ if and only if $\mathrm{Im} \, c \subseteq K$.  It follows
from this that the image of $\partial K$ under the map $\partial K \rightarrow \partial X_F$ is
precisely $\Delta_{L} \cup \Delta_{R} \cup \Delta_{L-R}$.  The second statement now follows from the fact that
$\partial K$ is isometric to $[0, \pi/2]$.

Now suppose that $c \not \in \Delta_{\infty}$ 
is a geodesic ray in $X_F$ issuing from $\ast$.  If $c$ is fixed by all of $F$ under the
action $\ast$, then the argument
of Section \ref{main} shows that $c \in \Delta_{L} \cup \Delta_{R} \cup \Delta_{L-R}$.

Conversely, suppose that 
$c \in \Delta_{L} \cup \Delta_{R} \cup \Delta_{L-R}$.  It follows from this and Lemma \ref{trav} that
$~Im~ c \subseteq K$, so every point $x \in ~Im~ c$ is within $1 + 2\sqrt{2}$ of a point in
$\{ 2, 3, \ldots \} \times \{ 2, 3, \ldots \} \subseteq \mathbb{R}^{2,+}$.  
We let $\widehat{T}_{m,n}$ denote
the tree in Figure \ref{done}, which consists of $T_{m,n}$ and one additional caret:

\begin{figure} [!h]
\begin{center}
\includegraphics{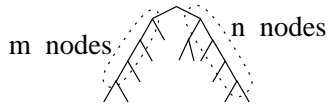}
\caption{A picture of the tree $\widehat{T}_{m,n}$, which consists of $T_{m,n}$ and one
additional caret.}   
\label{done}
\end{center}
\end{figure}

If $m,n \geq 2$, it is routine to check that 
\begin{eqnarray*}
x_0 \cdot T_{m,n} & = & T_{m+1, n-1}; \\
x_1 \cdot T_{m,n} & = & \widehat{T}_{m,n-1}.
\end{eqnarray*}
The trees $T_{m,n}$ and $T_{m+1, n-1}$ can be joined by an edge-path of length $2$ in $X_F$.  
The same goes for
$T_{m,n}$ and $\widehat{T}_{m,n-1}$, so $d( x_i \cdot T_{m,n} , T_{m,n} ) \leq 2$ $(i=0,1)$.

Let $t \geq 0$.  We have:
\begin{eqnarray*}
d( c(t) , x_i \cdot c(t) ) & \leq & d(c(t), T_{m,n}) + d(T_{m,n}, x_i \cdot T_{m,n} ) + 
d( x_i \cdot T_{m,n} , x_i \cdot c(t)) \\
& \leq & 4 + 4 \sqrt{2}.
\end{eqnarray*}

Since this estimate doesn't depend upon $t$, it follows that $x_0$, $x_1$ both fix $c$  under the action
$\ast$.  This implies that $F$ fixes $c$,
since $x_0$, $x_1$ generate $F$.
\end{proof}


\section{The Case of $\Delta_{\infty}$} \label{infty}


This section is devoted to an investigation of fixed points in
$\Delta_{\infty}$.  Our main result is the following:

\begin{theorem} \label{good}
The profile $\Delta_{\infty}$ contains no global fixed point of $T$ or $V$.
\end{theorem}

An immediate consequence of Theorem \ref{good} and Theorem \ref{fixTV} is:

\begin{corollary}
Thompson's groups $T$ and $V$ act without global fixed points on the boundaries-at-infinity of
their respective picture complexes.
\qed
\end{corollary}

The case for Thompson's group $F$ is more complicated.  I don't know whether 
$\Delta_{\infty} \subseteq \partial X_F$ contains fixed points of $F$ or not.
Example \ref{clue1} and Proposition \ref{clue2} give evidence for and against the existence
of fixed points in $\Delta_{\infty}$, respectively.

Our arguments will use a simple procedure for embedding any CAT(0) cubical complex into Hilbert space. 
\begin{proposition}
If $X$ is a CAT(0) cubical complex, then there is an embedding $\rho : X \rightarrow \ell^{2} \left( \mathcal{H} \right)$,
where $\mathcal{H}$ is the collection of hyperplanes in $X$.  The map $\rho$ doesn't increase distances.
\end{proposition}

\begin{proof}
Let $X$ be a CAT(0) cubical complex with a distinguished base vertex $\ast$.  For any given hyperplane $H$ in $X$, we
identify the closed $1/2$-neighborhood of $H$ with $H \times [0,1]$ in such a way that
$d( \ast, H \times \{ 0 \} ) < d ( \ast, H \times \{ 1 \} )$, and let $H_t$ denote $H \times \{ t \}$ for
$t \in [0,1]$.

Let $\rho : X \rightarrow \ell^{2} ( \mathcal{H} )$ send $x$ to
$\Sigma_{H \in \mathcal{H}} f_{H}(x) H$, where
$$f_{H}(x) = ~sup~ \left( \{ t \in [0,1] \mid [ \ast, x ] \cap H_t \neq \emptyset \} \cup \{ 0 \} \right).$$
It is clear that $\rho$ is an embedding, that each cube of $x$ is embedded into $\ell^{2}( \mathcal{H})$ isometrically,
and that the restrictions $\rho_{\mid C}$ (where $C$ is a cube) agree on overlaps.  Note also that each sum
$\Sigma_{H \in \mathcal{H}} f_{H}(x) H$ is finite, since any pair of vertices in $X$ are separated by at most
finitely many hyperplanes.  The first statement follows.

Let $[x,y]$ be a geodesic in $X$.  Since the restriction of $\rho$ to each cube is an isometry, $\rho$ preserves the
lengths of paths.  Therefore,
$$ d_{X}(x,y) = \ell \left( \rho [x,y] \right) \geq d_{\ell^{2}(\mathcal{H})} \left( \rho(x) , \rho(y) \right).$$
\end{proof}

\begin{example} \label{hilbert}
We consider the image of a certain $x \in X_F$ under the map $\rho$.
\begin{figure} [!h]  
\begin{center}
\includegraphics{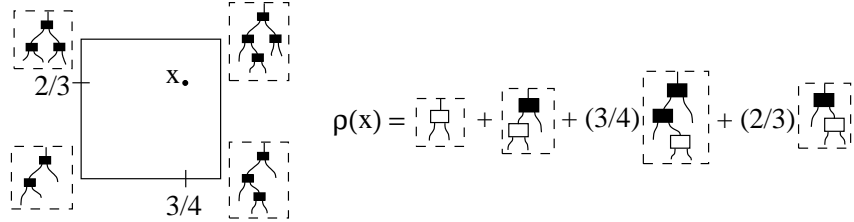}
\end{center}
\caption{A point $x$ in the cubical complex for Thompson's group $F$ (left) and its image under the map $\rho$ (right).}
\label{embed}
\end{figure}
Note that each hyperplane occurring in the sum $\rho(x)$ corresponds in a straightforward way
to a particular transistor $\Delta_{\infty}$.  We can use this fact to simplify our notation -- compare 
Figure \ref{simple} to the right
half of Figure \ref{embed}.
\begin{figure} [!h]
\begin{center}
\includegraphics{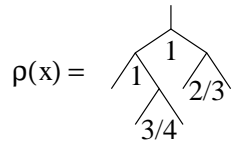}
\end{center}
\caption{A simpler notation for $\rho(x)$.}
\label{simple}
\end{figure}
In this way, we identify each point of $\mathrm{Im} \, \rho$ with a picture $\Delta$ such that
each maximal transistor is labelled by a number $t \in (0,1]$, and every other transistor
is labelled by the number $1$.  We will usually just omit the label of a transistor $T$ if $T$ is not maximal.  

With this convention, an element of $F$, $T$, or $V$ acts on a point of $\mathrm{Im} \, \rho$ by the usual picture multiplication.
This action is somewhat tricky to describe if a transistor from the acting element forms a dipole with an
exposed transistor labelled by a number $t \neq 1$.  In practice, however, we will always be able to avoid 
considering this
situation.  If there are no such dipoles, then the action is simple to describe: concatenate and reduce dipoles.
\end{example}

It will be helpful to have a vocabulary for describing subtrees of a given labelled tree $T$.  If $T_{\infty}$ is the full 
ordered rooted binary tree of infinite depth, we use binary strings to denote the vertices of degree three in $T_{\infty}$.
We give each edge in $T_{\infty}$ a label of $0$ or $1$; the label is $0$ if the edge forms the left half of a caret, and
$1$ if the edge forms the right half of a caret.  Now label each vertex $v$ in $T_{\infty}$ by the label of the unique geodesic
path from the root to $v$.  For instance, if the geodesic path from the root to $v$ passes through the right half of a caret twice,
and then through the left half of a caret, and then finally through the right half of a caret again, then the label of $v$ is $1101$.
The root has the empty label.  

Now if $T$ is an arbitrary labelled subtree (as in Figure \ref{simple}) 
of $T_{\infty}$ and $\mathrm{bin}$ is a binary string, we let $T_{\mathrm{bin}}$  
denote the labelled tree having the vertex $\mathrm{bin}$ as its root.  For instance, if $T$ is the labelled tree in Figure \ref{simple},
then $T_{1}$ consists of a single caret, labelled by the number $2/3$.  The tree $T_{01}$ is a single caret labelled by $3/4$.    

The following partial result will be used in the proof of Theorem \ref{good}.

\begin{proposition} \label{clue2}
Let $c: [0, \infty) \rightarrow X$ be a geodesic ray, where $X= \kpw$, $\kapw$, or $\kbpw$, 
$\mathcal{P} = \langle x \mid x = x^{2} \rangle$, and $w=x$.  Let $T$ be some rooted ordered binary tree in which each transistor is
labelled by a $1$ (and thus $T$ corresponds to a vertex in $\kpw$, $\kapw$, or $\kbpw$, as the case may be) such that the subtrees 
$T_{10}$ and $T_{11}$ each contain at least one caret.   

If $c(t) = T$ for some $t \in [0, \infty)$, then $x_{0} \ast c \neq c$.
\end{proposition}

\begin{proof}
We can express $\rho(c(t))$ as a tree of the form in Figure \ref{sec6arg} a), 
\begin{figure} [!h]
\begin{center}
\includegraphics{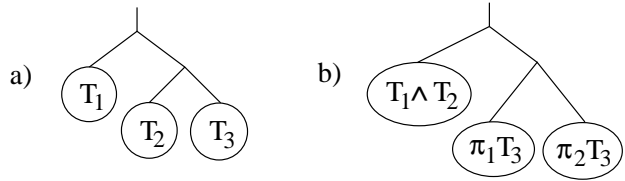}
\end{center}
\caption{a)  The form of $\rho(c(t))$.  The trees $T_2$ and $T_3$ both contain at least one caret.  b) the effect of letting
$x_0$ act on $\rho(c(t))$.}
\label{sec6arg}
\end{figure}
where each tree $T_{i}$, $i=1,2,3$, is a labelled tree, i.e., a picture in which every transistor is an 
$(x, x^{2})$-transistor, and the trees $T_2$ and $T_3$ each have at least one transistor.  
The effect of the action by $x_0$
is to transform the tree in Figure \ref{sec6arg}a) into the tree in Figure \ref{sec6arg}b).

Figure \ref{sec6arg}b) also serves as a definition of three different operations on labelled trees.  We now make this more explicit.
If $T$ is a labelled tree having at least one transistor labelled $1$, then the act of removing the topmost transistor
of $T$ leaves an ordered pair of trees $(\pi_1 T , \pi_2 T)$.  If $T_1$ and $T_2$ are labelled trees, then
$T_1 \wedge T_2$ is the unique labelled tree satisfying $\pi_1 ( T_1 \wedge T_2 ) = T_1$ and
$\pi_2 ( T_1 \wedge T_2 )$.

\begin{figure}
\begin{center}
\includegraphics{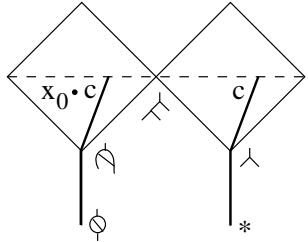}
\caption{On the right, we have the general picture describing a geodesic ray $c$ issuing from the basepoint $\ast$.  Every such
geodesic ray must cross the dotted horizontal line; at the moment it does so, it is precisely $\sqrt{2}$ units distant from 
its translate $x_{0} \cdot c$.}
\label{contra}
\end{center}
\end{figure}

Note that 
$$ || T - x_{0} \cdot T ||_{2}^{2} = || T_{1} - T_{1} \wedge T_{2} ||_{2}^{2} + || T_{2} - \pi_{1} T_{3} ||_{2}^{2} + ||T_{3} - \pi_{2} T_{3}||_{2}^{2}.$$
Since each caret in each tree is labelled with a $1$, each of the terms on the right side of the equation is an integer which counts the number of carets 
that are in one tree but not in the other.  It follows that
$|| T_{1} - T_{1} \wedge T_{2} ||_{2}^{2} \geq 2$ and $||T_{3} - \pi_{2}T_{3}||_{2}^{2} \geq 1$.  This implies that
$$d(c(t), x \cdot c(t)) \geq || T- x_{0} \cdot T ||_{2} \geq \sqrt{3}.$$

Now we appeal to Lemma \ref{monotone}(2), which implies that if there is $t' < t$ such that $d( c(t') , x_{0} \cdot c(t')) < \sqrt{3}$, then
$c$ and $x_{0} \cdot c$ are not asymptotic, i.e., $c \neq x_{0} \ast c$.  We produce a $t' < t$ 
where $d( c(t'), x_{0} \cdot c(t') ) = \sqrt{2}$
in Figure \ref{contra}. 

\end{proof}

\begin{remark}  \label{rem1}
\noindent (1)  It is reasonably clear from the proof that there are variations on this Proposition, in which the hypothesis that
$T_{10}$ and $T_{11}$ are non-trivial is replaced by similar assumptions on different subtrees.

\noindent (2)  
Let $v$ be a vertex in the cubical complex for Thompson's group $F$ such that
$\rho(v)$ is an ordered labelled rooted binary tree $\widehat{T}$.  Note that 
the coefficient of each caret in $\widehat{T}$ is either $1$ or $0$, since $v$ is a vertex.  
The above argument shows that
$$ \ast)  \quad d(v, x_{0} \cdot v) \geq \sqrt{ \left( 1 + ( \# \mathrm{~of~carets~in~} \widehat{T}_{10} ) \right) + 
\left( 1 + ( \# \mathrm{~of~carets~in~} \widehat{T}_{11} \right)},$$         
provided the trees $\widehat{T}_{10}$ and $\widehat{T}_{11}$ both contain at least one caret. 

This suggests a strategy for proving that $\Delta_{\infty}$ contains no fixed points of Thompson's group $F$.
Suppose that $c$ is a geodesic ray having the profile $\Delta_{\infty}$.  It follows from Lemma \ref{trav}
that $\rho(c(t))$ is a labelled ordered rooted binary tree $T(t)$, for any $t \geq 0$.  
For any given caret $C \in T_{\infty}$, $C$ occurs in $T(t)$ with the coefficient $1$
for $t$ sufficiently large, 
since $c$ crosses every hyperplane
in $\Delta_{\infty}$. 
If $c(t)$ passed close to a vertex $v$ for a large value of $t$, then the
inequality $\ast)$ shows that $d( v, x_{0} \cdot v)$ would be large, so that $d( c(t), x_{0} \cdot c(t))$ would
also be large, and thus $c \neq x_{0} \ast c$.  (Indeed, the proof of Proposition \ref{clue2} shows it
suffices to prove that $d( c(t), x_{0} \cdot c(t)) > \sqrt{2}$.)  Unfortunately, I know of no way to
control the distance of $c(t)$ from a vertex, since $c$ will generally travel through cubes of higher and higher
dimension, whose diameters go to infinity. 
\end{remark}  

\noindent \emph{Proof of Theorem \ref{good}:} 
Let $c \in \Delta_{\infty}$ be a geodesic ray in either $X_{T}$ or $X_{V}$
which represents a fixed point at infinity. 
We note first that $c(1)$ is the vertex labelled
by the finite tree $T_1$. If $\pi_1 \in T \subseteq V$ is as in Figure \ref{TV}b), it is clear that
$\pi_{1} \cdot T_1 = T_1$.  Now since $(\pi \cdot c)(1) = c(1)$ and 
$\pi \cdot c$, $c$ are asymptotic by our assumptions, 
it follows from Lemma \ref{monotone}(2) 
that $(\pi \cdot c)(t) = c(t)$, for $t \geq 1$.  

Now we refer to Figure \ref{contra}.  The element $\pi_1$ flips the diamond
on the right along its vertical axis, stabilizing the vertices $T_1$ and $T_2$.
Since the geodesic ray travels through this diamond and $(\pi_1 \cdot c)(t) = c(t)$
for $t \geq 1$, it must be that the ray $c$ travels along the straight line
connecting $T_1$ to $T_2$, which is vertical in the Figure.

Next, we consider the action of $\pi_2 \in T$ (see Figure \ref{wlog}b)) on $c$.  Note
first that $\pi_2$ stabilizes the vertex $T_2$.  It follows from this, Lemma \ref{monotone}(2), the 
equality $c( 1 + \sqrt{2}) = T_2$, and the assumption that $\pi_{2}$ fixes $c$, 
that $( \pi_{2} \cdot c )(t) = c(t)$ for $t \geq 1 + \sqrt{2}$.  
Lemma \ref{trav}
implies that $c( 1 + \sqrt{2} + \epsilon )$ is in the four dimensional cube $C$ having $T_2$ as its minimal
vertex and $T_3$ as its maximal vertex.  If we identify the cube $C$ with $[0,1]^{4}$ in such a way that
$(0,0,0,0) = T_2$ and $(1,1,1,1) = T_3$, then $\pi_2$ acts by cyclically permuting the coordinates of $C$.  It follows
from this that the geodesic ray travels along the diagonal of $C$ from $T_2$ to $T_3$.  In particular,
$c$ passes through $T_3$.

It now follows from Proposition \ref{clue2} that $c \neq x_0 \ast c$, which is a contradiction.
\qed 

\vspace{20pt}

We conclude with an example giving some evidence that
there may be fixed points of $F$ in $\Delta_{\infty}$. 

\begin{figure} [!h]
\begin{center}
\includegraphics{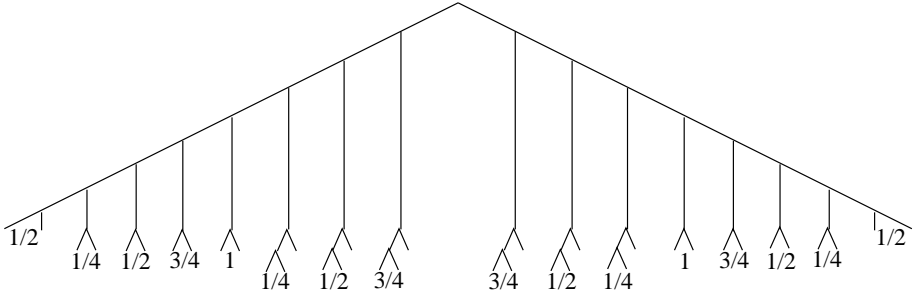}
\caption{This labelled tree represents a point in $X_F$ which is moved only a 
small distance by the generators of $F$. 
It should be possible to 
build a sequence of similar, arbitrarily large labelled trees which converge to a fixed point at infinity.}  
\label{trxm}
\end{center}  
\end{figure}

\begin{example} \label{clue1}       
Figure \ref{trxm} depicts a labelled tree $T$ (i.e., point in Hilbert space)
such that $||T - x_{0} \cdot T||_{2} = \sqrt{7}/2$.    The check is left as an exercise.
It is not difficult to see that $x_0$ acts on (most) of the trees along the bottom by a leftward shift.
In particular, the tree $T_{1\ldots 10}$ (where there are $n$ ones and a single $0$) is mapped
to $T_{1 \ldots 10}$ (where there are $n-1$ ones); the tree $T_{0 \ldots 01}$ ($n$ zeros) is mapped
to $T_{0 \ldots 01}$ ($n+1$ zeros).  
This suggests a principle for building larger trees that are moved only a small distance by $x_0$:
Begin with the tree $T_{m,m}$ and attach new trees $T'$ to the leaves, making sure that the tree attached 
at a given leaf is within $\epsilon$ of its neighbor to the immediate left, where $\epsilon$ will depend on
$m$.  The object is to make sure that $|| T - x_{0} \cdot T ||_{2}$ is less than or equal to $\sqrt{2}$, and
it is not too difficult to see that this can be done for any $m$.  (A little experimentation also shows that
it is useful to label the leftmost and rightmost carets of $T_{m,m}$ with $1/2$.)  Moreover, the trees $T'$
that are attached ``close" to the root (in a sense that depends on $m$)
can be made arbitrarily large.  Thus, we can make a sequence of labelled
trees $\overline{T}_{k}$, which are each moved less than $\sqrt{2}$ units in the Hilbert metric, and 
gradually fill up the complete binary tree $T_{\infty}$.  One can then hope  that the corresponding
points $z_k$ in $X_F$ are also moved only a small distance by $x_0$ (as seems likely), so that  some subsequence of $z_k$
converges to a point $\zeta$ at infinity which is fixed by the action of $x_0$.  With additional care, it should also be
possible to do this so that each point $z_k$ is likewise moved only a small distance by $x_1$, and therefore $\zeta$ would
be fixed by all of $F$.     

It seems very likely that all of the above can be done.  This is not enough, however, because the point $\zeta$ may well
fail to have the profile $\Delta_{\infty}$.  Indeed, the tree in Figure \ref{trxm} has the property that most of its norm
is contributed by $T_{7,7}$.  It appears likely that any tree in the 
sequence $\overline{T}_k$ will have most of its norm contributed by $T_{n,n}$ (for appropriate $n$), and this may mean
that $\zeta \in \Delta_{L} \cup \Delta_{L-R} \cup \Delta_{R}$.  I conjecture the following:

\begin{conjecture}
A point $\zeta \in \partial X_F$ is fixed by all of $F$ if and only if $\zeta \in \Delta_{L} \cup \Delta_{L-R} \cup \Delta_{R}$.
In particular, the only fixed points for the action of $F$ on $\partial X_F$ lie on an arc of Tits length $\pi /2$. 
\end{conjecture}
\end{example}


\bibliography{infty}
\nocite{*}
\bibliographystyle{plain}

\end{document}